\documentclass[a4paper,11pt]{amsart}

\usepackage{geometry}
\usepackage{amsmath}
\usepackage{amsfonts}
\usepackage{amssymb}
\usepackage{amsthm}
\usepackage{mathrsfs}
\usepackage{latexsym}

\usepackage[english]{babel}

\usepackage{xfrac}

\usepackage{mathtools}

\usepackage[latin1]{inputenc}
\usepackage{graphicx}

\usepackage[pdftex]{hyperref}

\usepackage{xypic}
\input xy
\xyoption{all}

\theoremstyle{plain}
\newtheorem{thm}{Theorem}
\newtheorem{lem}[thm]{Lemma}
\newtheorem{prop}[thm]{Proposition} 
\newtheorem{cor}[thm]{Corollary} 

\theoremstyle{definition} 
\newtheorem{defi}[thm]{Definition} 

\theoremstyle{remark} 
\newtheorem{remark}[thm]{Remark}

\title[Heisenberg algebra and framed sheaves]{Representations of the Heisenberg algebra and moduli spaces of framed sheaves}

\author{Francesco Sala}

\address{Mathematical Physics Sector, SISSA, Via Bonomea 265, 34100, Trieste (ITALY)\newline
Laboratoire Paul Painlevé, Université Lille 1, Cité Scientifique, 59655, Villeneuve D'Ascq Cedex (FRANCE)}

\email{francesco.sala@sissa.it}

\author{Pietro Tortella}

\address{Mathematical Physics Sector, SISSA, Via Bonomea 265, 34100, Trieste (ITALY)\newline
Laboratoire Paul Painlevé, Université Lille 1, Cité Scientifique, 59655, Villeneuve D'Ascq Cedex (FRANCE)}

\email{pietro.tortella@sissa.it}

\begin{document}

\begin{abstract}
Let $\mathcal{M}_{X,D}(r,c,n)$ be the moduli spaces of torsion free sheaves on a complex smooth connected projective surface $X$, framed along a smooth connected genus zero curve $D.$ This paper gives a `geometrical' construction of an action of a Heisenberg algebra on the homology of $\mathcal{M}_{X,D}(r,c,n)$ (more precisely, on the direct sum $\bigoplus_n \mathrm{H}_* (\mathcal{M}_{X,D}(r,c,n))$) using correspondences. This result generalizes Nakajima's construction for the Hilbert schemes of points.
\end{abstract}

\maketitle

\section{Introduction}

Let $X$ be a complex smooth connected quasi-projective surface and $X^{[n]}$ the Hilbert scheme of length $n$ zero-dimensional subschemes of $X$. In \cite{book:naka} Nakajima builds a representation of the Heisenberg-Clifford super-algebra on the direct sum of the homology groups $\bigoplus_{n\geq 0} \mathrm{H}_*\left(X^{[n]}\right)$, showing that the generating function of the Poincaré polynomials, found by G\"{o}ttsche and Soergel in \cite{art:gott2}, coincides with the character formula of the representation. Nakajima's result gives a supporting evidence to the `S-duality conjecture' proposed by Vafa and Witten in \cite{art:vawi}.

The Hilbert scheme $X^{[n]}$ can be viewed as the moduli space of rank one sheaves over $X$ with trivial determinant and second Chern class equals to $n$. Let $X$ be a smooth projective surface; in \cite{art:bara2} Baranovsky generalizes Nakajima's result to the moduli spaces $\mathcal{M}_X(r,L,n)$ of Gieseker-stable torsion free sheaves on $X$ of rank $r$, determinant $L$ and second Chern class $n.$

Let $D$ be a big and nef curve in $X$. In \cite{art:moduli_framed} Bruzzo and Markushevich built a fine moduli space $\mathcal{M}_{X,D}(r,c,n)$ for torsion free sheves on $X$ with invariants $(r,c,n)$, framed along the divisor $D.$ The Hilbert scheme of points $(X\setminus D)^{[n]}$ can be viewed as the rank one case of this moduli space.

Moreover moduli spaces of framed sheaves are studied because they provide a desingularization of the moduli spaces of ideal instantons, so their equivariant cohomology under suitable toric actions is relevant to the computation of partition functions in topological quantum field theory (e.g. see \cite{art:nekra}, \cite{art:brutanfu}, \cite{art:nakayo} and \cite{art:brutan}).

The aim of this work is to generalize Nakajima's construction to these moduli spaces: we explicitly construct \emph{Nakajima's operators} on the direct sum of the homology groups $\bigoplus_{n} \mathrm{H}_* (\mathcal{M}_{X,D}(r,c,n))$ and show that they satisfy the Heisenberg commutation relations and hence we get a representation of the Heisenberg algebra generated by \emph{Nakajima's operators} on this space.

Using this construction one can expect to generalize to higher rank case Carlsson and Okunkov result (\cite{art:carl}): they use Nakajima's operators for Hilbert schemes $X^{[n]}$ to define a class of vertex operators, whose trace gives the Nekrasov's partition function in the case we have a toric action on $X.$ This could be relevant to provide a mathematical motivation for the conjecture about a correspondence between the Nekrasov's partition function for supersymmetric Yang-Mills topological field theories and the conformal blocks of a 2-dimensional conformal field theory (e.g. see \cite{art:agata}).

The article is structured in the following way: in section 2 we introduce framed sheaves over a smooth connected projective surface and we outline how one can construct moduli spaces of these; in section 3 we construct the \textit{Nakajima operators}, wich will lead us to the representation we are seeking; in section 4 we state and prove a generalization of Nakajima's theorem for the case of moduli spaces of framed sheaves, postponing to the last two sections a proof of the more technical propositions.

\subsection*{Technical remark} 

In this paper we shall assume that the framing divisor $D$ is a smooth connected genus zero curve and this forces $X$ to be a rational surface; we assume this hypothesis because we need the moduli spaces of framed sheaves to be smooth varieties and for the technical lemma \ref{lem:arrow}. A priori, one can replace this assumption by `\emph{let} $X$ \emph{and} $D$ \emph{be such that the moduli spaces} $\mathcal{M}_{X,D}(r,c,n)$ \emph{are smooth}' and one can construct the operators in the same way, but we do not know how to prove the Heisenberg commutator relations in this general case. Anyway the main explicit examples of framed sheaves and their moduli (e.g. see chapter 2 in \cite{book:naka}, \cite{art:henni}, \cite{phd:rava}) satisfy our hypothesis.

\subsection*{Conventions and notation} 

All schemes we are dealing with are of finite type over $\mathbb{C}$, by `variety' we mean a reduced separated scheme; a `sheaf' is always coherent.

\subsection*{Acknowledgements}

We would like to thank Ugo Bruzzo for suggesting us this problem and for constant support. Thanks to Vladimir Baranovsky for useful remarks about his article and suggestions, Claudio Rava and Emanuele Macrì for useful conversations.
	
\section{Framed Sheaves and their Moduli}	

In this section we introduce the notion of framed sheaf and we give a construction of moduli spaces of these objects.

\subsection*{Generalities}

Let $X$ be a smooth connected projective surface. Fix a smooth connected curve $D$ in $X$.

\begin{defi}
A \emph{framed sheaf of rank} $r$ on $X$ with framing divisor $D$ is a pair $\mathcal{E}=(E,\phi^{E})$ where $E$ is a torsion free sheaf of rank $r$ on $X$, $\phi^{E}$ is a morphism $\phi^{E}:E\rightarrow O_D^{\oplus r}$ such that $\phi^{E}\vert_{D}$ is an isomorphism. If $(E,\phi^{E})$ and $(F,\phi^{F})$ are two framed sheaves of the same rank $r$, we define a \emph{morphism of framed sheaves} as a morphism $\alpha:E\rightarrow F$ such that for some $\lambda\in\mathbb{C}$ we have a commutative diagram
\begin{equation}
 \xymatrix{	E\ar^{\alpha}[r]\ar^{\phi^{E}}[d] & F\ar^{\phi^{F}}[d] \\
		O_D^{\oplus r} \ar^{\cdot \lambda}[r] & O_D^{\oplus r}
	}
\end{equation}
\end{defi}

Clearly the set of morphisms of framed sheaves is a linear subspace of the space of morphisms between the corresponding underlying sheaves.

Now let $(E,\phi^{E})$ be a framed sheaf of rank $r$ on $X.$ Because of the framing, the stalks at the points on the divisor $D$ are free of rank $r$, so $E$ is locally free in a neighborhood of $D$, hence in this neighborhood the sheaves $E$ and $E^{\vee\vee}$ are isomorphic. Thus we have a natural framing for $E^{\vee\vee}$. Moreover the support of $A$ is disjoint from $D$. In the following we denote by $\mathcal{E}^{\vee\vee}$ the framed sheaf $(E^{\vee\vee},\phi^{E})$. Note that the inclusion of $E$ in $E^{\vee\vee}$ induces a morphism between the corresponding framed sheaves.

\subsection*{Construction of the moduli space}

We now come to the construction of the moduli space of framed sheaves on a surface.

\begin{defi}
Let $S$ be a scheme. A \emph{family of rank} $r$ \emph{framed sheaves parametrized by} $S$ is a pair $(G,\phi^{G})$ where $G$ is a torsion free sheaf of rank $r$ on $X\times S$ flat over $S$, $\phi^{G}$ is a morphism $\phi^{G}:G\rightarrow p_X^* O_D^{\oplus r}$ such that  $\phi^{G}\vert_{D}:G\vert_{D\times S}\rightarrow p_X^* O_D^{\oplus r}$ is an isomorphism.

Remark that a family of rank $r$ framed sheaves over $X$ parametrized by $S$ with framing divisor $D$ is in particular a rank $r$ framed sheaf over the product $X\times S$ with framing divisor $D\times S.$

We say that two families $\mathcal{G}$ and $\mathcal{G}'$ of rank $r$ framed sheaves over $X$ parametrized by $S$ are isomorphic if they are isomorphic as rank $r$ framed sheaves over $X\times S$.
\end{defi}

Now fix invariants $r,n\in\mathbb{Z}$ with $r>0$ and $c\in \mathrm{NS}(X)$; we can define the functor
\begin{equation*}
\underline{\mathcal{M}_{X,D}}(r,c,n) : (Schemes)^\circ \longrightarrow (Sets)
\end{equation*} 
that associates to every scheme $S$ the set of isomorphism classes of families $(G,\phi^{G})$ of rank $r$ framed sheaves parametrized by $S$ such that the fibres $G_s$ have rank $r$ and Chern classes $c$ and $n$, while to every morphism of schemes $f:T\rightarrow S$ associates the pull-back $f^*$ that sends families parametrized by $S$ in families parametrized by $T$.

\begin{remark}
A necessary condition for this functor to be nonempty is $\displaystyle\int_D c =0 $, indeed in general $\displaystyle \int_D c_1(E)= c_1(E\vert_{D})$ and if $E$ is trivial along $D$, then $c_1(E\vert_{D})=c_1(O_D^{\oplus r})=0$.
\end{remark}

In \cite{art:moduli_framed} Bruzzo and Markushevich give (in a more general setting) a proof of the representability of the functor $\underline{\mathcal{M}_{X,D}}(r,c,n)$; using their result we have the following:
\begin{thm}
Let $X$ be a smooth connected projective surface and $D$ a smooth connected curve in $X$ which is a big and nef divisor. Then there exists a quasi-projective separated scheme $\mathcal{M}_{X,D}(r,c,n)$ that represents the functor $\underline{\mathcal{M}_{X,D}}(r,c,n)$. 
\end{thm}

Note that $\mathcal{M}_{X,D}(r,c,n)$ is a \emph{fine} moduli space for framed sheaves on $X$ with framing divisor $D$ and topological invariants $r, c$ and $n$. The adjective `fine' means the existence of a \emph{universal} framed sheaf $\bar{\mathcal{E}}=(\bar{E}, \phi^{\bar{E}})$, i.e. a family of framed sheaves parametrized by $\mathcal{M}_{X,D}(r,c,n)$ with the following \emph{universal} property: for any family $(G, \phi^G)$ of framed sheaves parametrized by $S$ there exists a unique morphism $g:S\rightarrow \mathcal{M}_{X,D}(r,c,n)$ and an isomorphism $\alpha: G\rightarrow (\mathrm{id}\times g)^* (\bar{E})$ such that $\phi^G\vert_{D\times S}=((\mathrm{id}\times g)^*\phi^{\bar{E}}\circ \alpha)\vert_{D\times S}.$

From now on we suppose that $D$ is a smooth connected curve in $X$ that satisfies the following conditions:
\begin{itemize}
\item[(a)] $D$ is a big and nef divisor,
\item[(b)] $D\cong \mathbb{CP}^1.$
\end{itemize}

\begin{remark}\label{rem:ratio}
Since $D$ is isomorphic to the complex projective line and it has positive self intersection, $X$ is a rational connected surface, hence it is a rational surface. Moreover the arithmetic genera of $D$ and $X$ are zero and $\mathrm{H}_i(X,\mathbb{C})=0$ for $i=1,3.$
\end{remark}

From these conditions follows the lemma:
\begin{lem}\label{lem:morph_framed}
Let $(E,\phi^E)$ and $(F,\phi^F)$ be framed sheaves on $X$. Then
\begin{equation*}
\mathrm{Hom}(E,F)\cong \mathrm{End}(\mathbb{C}^r).
\end{equation*}
\end{lem}
\proof
Let $f$ be a non trivial morphism in $\mathrm{Hom}(E,F)$; since $E\vert_D$ and $F\vert_D$ are trivial, $f\vert_D$ is a constant matrix, i.e. is an element of $\mathrm{End(\mathbb{C}^r)}.$
 
Consider the complete linear system $\vert D\vert$ and let $Y$ be its base locus: we can define a morphism $\pi:X\setminus Y\rightarrow |D|$ that associates to each point $x$ of $X\setminus Y$ the hyperplane in $\vert D\vert$ consisting of those sections vanishing at $x.$ 

Let $V$ be the open set of divisor $C\in \vert D\vert$ such that $C$ is a smooth curve of genus zero. Let $U$ be the set consisting on curves $C\in V$ such that $E\vert_C$ and $F\vert_C$ are the trivial bundle. By lemma 2.3.1 in \cite{phd:king} and semicontinuity theorem $U$ is open. We have that $f\vert_C$ is a constant matrix for all $C\in U.$
So we get that $f$ is a constant matrix in $\pi^{-1}(U)$, that is a dense open subset of $X$. Thus $f$ is constant in the whole $X$ and the thesis is proved.
\endproof

Now we want to study the local properties of the moduli space, in particular we want a description of its tangent space. We begin with a vanishing lemma:
\begin{lem} \label{lem:vanishing_framed}
 Let $(E,\phi^{E})$ be a framed sheaf on $X$ with framing divisor $D$. Then 
\begin{itemize}
 \item[(a)] $\mathrm{Hom}(E,E(-D))=0$,
 \item[(b)] $\mathrm{Ext}^2(E,E(-D))=0.$
\end{itemize}
\end{lem}
\proof
 \textbf{(a)}  Let $f$ be a non trivial morphism in $\mathrm{Hom}(E,E(-D))$; since $E\vert_D$ is trivial we have that 
\begin{equation*}
 \mathrm{Hom}(E\vert_D,E(-D)\vert_D) \cong \mathrm{H}^0(D,\mathcal{H}om(O_D^{\oplus r},O_D^{\oplus r}\otimes O_X(-D))) \cong	\mathrm{H}^0(D,O_D(-a))^{\oplus r}
\end{equation*}
where $a=D.D$. Since $D$ is a big and nef divisor, we get that $a$ is positive, hence $\mathrm{H}^0(D,O_D(-a))=0$ and therefore $f$ is zero along $D$. 
 
As before, consider the morphism $\pi:X\setminus Y\rightarrow |D|.$ Let $T$ be the set of divisor $C\in \vert D\vert$ such that $E\vert_C$ is a semistable sheaf on $C$. Note that $T\neq \emptyset$ because $D$ belongs to $T$ and $T$ is open in $\vert D \vert$ because semistability is an open property. Since $p(E\vert_C, m)>p(E(-D)\vert_C, m)$ for all $C\in \vert D \vert$, by proposition 1.2.7 in \cite{book:hl} we get that $\mathrm{Hom}(E\vert_C, E(-D)\vert_C)=0$ for all $C\in T$ and therefore $f\vert_C=0$. So we can use the same argument of lemma \ref{lem:morph_framed} and therefore $f$ is zero in the whole $X.$

\textbf{(b)} Using Serre's duality we have 
\begin{equation*}
\mathrm{Ext}^2(E,E(-D)) \cong \mathrm{Hom}(E,E\otimes O_X(D)\otimes \omega_X)^\vee.
\end{equation*}
We have $\omega_D\cong (\omega_X\otimes O_X(D))\vert_D$ and $\mathrm{deg}\, \omega_D < 0$ by adjunction formula, so any element $f$ in $\mathrm{Hom}(E,E\otimes O_X(D)\otimes \omega_X)$ is zero when restricted to $D$; in the same way, for all $C\in T$ we get $\mathrm{deg}\, O_X(K_X+D)\vert_C < 0$ and using the same argument as before, the thesis follows.
\endproof

Thanks to this lemma and to theorem 4.3 in \cite{art:moduli_framed}, we get
\begin{thm}
 Let $\mathcal{E}=(E,\phi^{E})$ be a framed sheaf with invariants $r,c$ and $n$. Then the Zariski tangent space at the corresponding point $[\mathcal{E}]$ of $\mathcal{M}_{X,D}(r,c,n)$ is
\begin{equation*}
 T_{[\mathcal{E}]}\mathcal{M}_{X,D}(r,c,n) \cong \mathrm{Ext}^1(E,E(-D)).
\end{equation*}
Moreover $\mathcal{M}_{X,D}(r,c,n)$ is a smooth quasi-projective variety of dimension
\begin{equation*}
 \dim_{\mathbb{C}} \mathcal{M}_{X,D}(r,c,n) = 2rn - (r-1) \int_X c^2.
\end{equation*}
\end{thm}

From now on we put $\displaystyle b:=- (r-1) \int_X c^2.$

\begin{remark}
Note that there does not exist a general theorem about the nonemptiness of these moduli spaces. There are few known cases of nonempty moduli spaces of framed sheaves. The simplest case is given by $r=1$, $c=0$ and $n>0$: the moduli space $\mathcal{M}_{X,D}(1,0,n)$ is isomorphic to the Hilbert scheme of points $X_0^{[n]}$, where $X_0=X\setminus D$. This space is a smooth quasi-projective variety of dimension $2n$. Another case is given by $X=\mathbb{P}^2$ and $D$ a line in $\mathbb{P}^2$ (see for example \cite{book:naka}, chapter II); in this case one shows that $\mathcal{M}_{\mathbb{P}^2,D}(r,c,n)$ is nonempty if and only if $c=0$ and $n>0$. A recently known case is one on the Hirzebruch surfaces: let $X=\mathbb{F}_p$ be the $p$-th Hirzebruch surface and $D$ be the `line at infinity', let $c= aE$, where $a\in\mathbb{Z}$, $0\leq a<r-1$, and $E$ the only curve in $\mathbb{F}_p$ with negative self-intersection, in \cite{phd:rava} Rava proves that the moduli space $\mathcal{M}_{\mathbb{F}_p,D}(r,c,n)$ is nonempty if and only if the number $n+\frac{r-1}{2r}pa^2$ is an integer and $n+\frac{1}{2}pa(a-1)\geq 0$. Moreover he proves that when equality holds, $\mathcal{M}_{\mathbb{F}_1,D}(r,c,n)=Grass(a,r)$ and $\mathcal{M}_{\mathbb{F}_p,D}(r,c,n)=T(Grass(a,r))^{\oplus (p-1)}$ for $p>1$, where $Grass(a,r)$ is the Grassmannian of $a$-planes inside $\mathbb{C}^r.$ 
\end{remark}

In the following we assume that $X,D,r,c,n$ are such that the moduli space $\mathcal{M}_{X,D}(r,c,n)$ is nonempty.

\subsection*{Tangent bundle of $\mathcal{M}_{X,D}(r,c,n)$}

First of all we want to define the \textit{Kodaira-Spencer map} for framed sheaves in terms of cocycles. Let S be a scheme. 

Let $p_S$ and $p_X$ be the projections from $Y=X\times S$ to $S$ and $X$ respectively. Consider a family $(G, \phi^G)$ of framed sheaves parametrized by $S$. Let $G^{\bullet}$ be a finite locally free resolution $G^{\bullet}\rightarrow G$.

Following chapter 10 in \cite{book:hl}, we can define the \textit{Atiyah class} $A(G^{\bullet})$ in terms of connections: let $p_1, p_2: Y\times Y\rightarrow Y$ be the projections to the two factors. Let $\mathcal{I}$ be the ideal sheaf of the diagonal $\Delta\subset Y\times Y$ and let $O_{2\Delta}$ denote the structure sheaf of the first infinitesimal neighbourhood of $\Delta$.

Choose an open affine covering $\,\mathcal{U}=\{U_i\vert i\in I\}$ such that the restriction of the sequence
\begin{equation*}
0\longrightarrow G^q\otimes \Omega_Y\longrightarrow (p_1)_{*}(p_2^{*}G^q\otimes O_{2\Delta})\longrightarrow G^q\longrightarrow 0
\end{equation*}
to $U_i$ splits for all $q$ and $i$. Thus there are local connections $\nabla_i^q: G^q\vert_{U_i}\rightarrow G^q\otimes \Omega_Y\vert_{U_i}$. Note that the difference of two local connections is an $O_Y$-linear map. Define cochains $\alpha'\in C^1(\mathcal{H}om^0(G^{\bullet}, G^{\bullet}\otimes \Omega_Y), \mathcal{U})$ and $\alpha''\in C^0(\mathcal{H}om^1(G^{\bullet}, G^{\bullet}\otimes \Omega_Y), \mathcal{U})$ as follows:
\begin{eqnarray*}
\alpha'^q_{i_0\,i_1}&=& \nabla^q_{i_0}\vert_{U_{i_0\,i_1}}-\nabla^q_{i_1}\vert_{U_{i_0\,i_1}},\\
\alpha''^q_i&=&\mathrm{d}_{G^{\bullet}}\circ \nabla_i^q-\nabla^{q+1}_i\circ \mathrm{d}_{G^{\bullet}},
\end{eqnarray*}
where $\mathrm{d}_{G^{\bullet}}$ is the differential of the complex $G^{\bullet}$. The element $\alpha=\alpha'+\alpha''$ is a cocycle in the total complex associated to the double complex $C^{\bullet}(\mathcal{H}om^{\bullet}(G^{\bullet}, G^{\bullet}\otimes \Omega_Y), \mathcal{U})$. The cohomology class of $\alpha$ in $\mathbb{E}\mathrm{xt}^1(G^{\bullet}, G^{\bullet}\otimes \Omega_Y)$ is the \textit{Atiyah class} $A(G^{\bullet})$ of the sheaf $G$.

Now, we want to care about the morphism $\phi^G$. The \emph{connecting} morphism $\varepsilon: G^0\rightarrow G$ induces a cohomology class $[\varepsilon\otimes \mathrm{id}_{\Omega_Y}]$ in $\mathbb{E}\mathrm{xt}^0(G^{\bullet}\otimes \Omega_Y, G\otimes \Omega_Y)$ and therefore $[\overline{\alpha}]:=A(G^{\bullet})\otimes [\varepsilon\otimes \mathrm{id}_{\Omega_Y}]$ is a cohomology class in $\mathbb{E}\mathrm{xt}^1(G^{\bullet}, G\otimes \Omega_Y)$, where $\otimes$ is the Yoneda product for Ext-groups of complexes of sheaves.

Since $\mathcal{H}om^0(G^{\bullet}, G\otimes \Omega_Y)=\mathcal{H}om(G^0, G\otimes \Omega_Y)$ and $\mathcal{H}om^n(G^{\bullet}, G\otimes \Omega_Y)=0$ for $n>0$, we have that $\overline{\alpha}$ is an element in
\begin{equation*}
C^1(\mathcal{H}om(G^0, G\otimes \Omega_Y))\subset C^1(\mathcal{H}om^{0}(G^{\bullet}, G\otimes \Omega_Y\stackrel{\phi^G\otimes \mathrm{id}}{\longrightarrow}p_X^{*}O_D^{\oplus r}\otimes \Omega_Y)).
\end{equation*}

We define morphisms $\nabla_i: G^0\vert_{U_i}\rightarrow (G\otimes \Omega_Y)\vert_{U_i}$ as compositions
\begin{equation*}
G^0\vert_{U_i} \stackrel{\nabla_i^0}{\longrightarrow}(G^0\otimes \Omega_Y)\vert_{U_i}\stackrel{\varepsilon\otimes id}{\longrightarrow} (G\otimes \Omega_Y)\vert_{U_i}
\end{equation*}
and morphisms $\widetilde{\phi_i}: G^0\vert_{U_i}\rightarrow (p_X^{*}O_D^{\oplus r}\otimes \Omega_Y)\vert_{U_i}$ as compositions
\begin{equation*}
G^0\vert_{U_i} \stackrel{\nabla_i}{\longrightarrow}(G\otimes \Omega_Y)\vert_{U_i} \stackrel{\phi^G\otimes \mathrm{id}}{\longrightarrow}(p_X^{*}O_D^{\oplus r}\otimes \Omega_Y)\vert_{U_i}
\end{equation*}
hence we obtain a cochain $\widetilde{\phi}\in C^0(\mathcal{H}om^1(G^{\bullet}, G\otimes \Omega_Y\stackrel{\phi^G\otimes \mathrm{id}}{\longrightarrow}p_X^{*}O_D^{\oplus r}\otimes \Omega_Y))$.

Since $\widetilde{\phi_{i_0}}\vert_{U_{i_0\,i_1}}-\widetilde{\phi_{i_1}}\vert_{U_{i_0\,i_1}}=((\phi^G\circ \varepsilon)\otimes \mathrm{id}_{\Omega_Y})\vert_{U_{i_0\,i_1}}\circ \alpha'^{\,0}_{i_0\,i_1}$, we have that $A(G,\phi^G):=\overline{\alpha}+\widetilde{\phi}$ defines an element of 
\begin{equation*}
\mathbb{E}\mathrm{xt}^1(G^{\bullet}, G\otimes \Omega_Y\stackrel{\phi^G\otimes \mathrm{id}}{\longrightarrow}p_X^{*}O_D^{\oplus r}\otimes \Omega_Y)=\mathbb{E}\mathrm{xt}^1(G, G\otimes \Omega_Y\stackrel{\phi^G\otimes \mathrm{id}}{\longrightarrow}p_X^{*}O_D^{\oplus r}\otimes \Omega_Y).
\end{equation*}

If we consider the induced section of $A(G,\phi^G)$ under the global-to-local map 
\begin{equation*}
\mathbb{E}\mathrm{xt}^1(G, G\otimes \Omega_Y\stackrel{\phi^G\otimes \mathrm{id}}{\longrightarrow} p_X^{*}O_D^{\oplus r}\otimes \Omega_Y)\rightarrow \mathrm{H}^0(S, \mathcal{E}xt^1_{p_S}(G, G\otimes \Omega_Y\stackrel{\phi^G\otimes \mathrm{id}}{\longrightarrow}p_X^{*}O_D^{\oplus r}\otimes \Omega_Y))
\end{equation*}
where $\mathcal{E}xt_{p_S}^1$ is the first right derived functor of $(p_S)_{*}\circ \mathcal{H}om$, then the direct sum decomposition $\Omega_Y=p_S^{*}\Omega_S\oplus p_X^{*}\Omega_X$ leads to an analogous decomposition $A(G,\phi^G)=A(G,\phi^G)'+A(G,\phi^G)''$. The \textit{Kodaira-Spencer map} associated to $(G, \phi^G)$ is
\begin{equation*}
KS: \Omega_S^\vee\stackrel{A(G,\phi^G)'}{\longrightarrow} \Omega_S^\vee\otimes\mathcal{E}xt^1_{p_S}(G, G\otimes p_S^{*}\Omega_S\stackrel{\phi^G\otimes\mathrm{id}}{\longrightarrow}p_X^{*}O_D^{\oplus r}\otimes p_S^*\Omega_S)\rightarrow \mathcal{E}xt^1_{p_S}(G, G\stackrel{\phi^G}{\rightarrow}p_X^{*}O_D^{\oplus r})
\end{equation*}

\begin{remark}
It is easy to prove that the complex $G\stackrel{\phi^G}{\rightarrow}p_X^{*}O_D^{\oplus r}$ is quasi-isomorphic to $G\otimes p_X^{*}O_X(-D)$, hence $\mathcal{E}xt^1_{p_S}(G, G\stackrel{\phi^G}{\rightarrow}p_X^{*}O_D^{\oplus r})\cong \mathcal{E}xt^1_{p_S}(G, G\otimes p_X^{*}O_X(-D))$.
\end{remark}

Now we are able to give the following characterization of the tangent bundle of $\mathcal{M}_{X,D}(r,c,n)$.

\begin{thm}\label{thm:tanbundle}
Let $(\bar{E},\phi^{\bar{E}})$ be the universal family of $\mathcal{M}_{X,D}(r,c,n)$ and let $p$ be the projection from $\mathcal{M}_{X,D}(r,c,n)\times X$ to $\mathcal{M}_{X,D}(r,c,n)$. The Kodaira-Spencer map for $\mathcal{M}_{X,D}(r,c,n)$ induces a canonical isomorphism
\begin{equation*}
\mathrm{T}_{\mathcal{M}_{X,D}(r,c,n)}\cong \mathcal{E}xt_p^1(\bar{E}, \bar{E}\otimes p_X^* O(-D)).
\end{equation*}
\end{thm}
\proof
To prove this theorem, we can use a similar argument of theorem 10.2.1 in \cite{book:hl}: by lemma \ref{lem:vanishing_framed} we get that the sheaf $\mathcal{E}xt_p^1(\bar{E}, \bar{E}\otimes p_X^* O(-D))$ commutes with base change and $\mathcal{E}xt_p^1(\bar{E}, \bar{E}\otimes p_X^* O(-D))$ is locally free.

Moreover, by example 10.1.9 in \cite{book:hl} we get that the Kodaira-Spencer map is an isomorphism on the fibres.
\endproof

\subsection*{Donaldson-Uhlenbeck partial compactification}

Let $\mathcal{M}^{reg}_{X,D}(r,c,n)$ be the open subset in $\mathcal{M}_{X,D}(r,c,n)$ consisting of framed sheaves $(F,\phi^{F})$ with $F$ locally free.

Now remember that each framed sheaf $\mathcal{E}=(E,\phi^{E})$ has a natural inclusion in $\mathcal{E}^{\vee\vee}=(E^{\vee\vee},\phi^{E})$ and the quotient $\sfrac{E^{\vee\vee}}{E}$ is a zero-dimensional sheaf whose support is disjoint from $D$. In this way to any framed sheaf $\mathcal{E}$ with invariants $r,c,n$ we can associate the pair $(\mathcal{E}^{\vee\vee},\sum_i a_i [x_i])$ where $x_i$ are the points of the support of the quotient and $a_i$ is the length of the quotient at $x_i$; remark that if $c_2(E^{\vee\vee})=s$ then $\mathrm{length}(\sfrac{E^{\vee\vee}}{E})=n-s$, hence $\sum_i a_i x_i$ is an element of the $(n-s)$-th symmetric product $S^{n-s} X$ and $\mathcal{E}^{\vee\vee}$ is an element of $\mathcal{M}^{reg}_{X,D}(r,c,s)$.

This association gives a map
\begin{equation}\label{eq:hilch}
\pi_r: \mathcal{M}_{X,D}(r,c,n)\longrightarrow \widehat{\mathcal{M}_{X,D}(r,c,n)}:=\coprod_{s=0}^n \mathcal{M}^{reg}_{X,D}(r,c,s)\times S^{n-s}X ;
\end{equation}
in \cite{art:uhlenbeck_framed} the right hand side space is given a structure of scheme (similar to the Donaldson-Uhlenbeck partial compactification of instantons) such that the map $\pi_r$ is a projective morphism. 

Now we are going to study the fibers of this map, that will be useful later on: let $\mathfrak{p}=((G,\phi^{G}),\sum_j m_j [x_j])$ be a point in $\mathcal{M}^{reg}_{X,D}(r,c,s)\times S^{n-s}X$, the fiber $\pi_r^{-1}(\mathfrak{p})$ parametrizes framed sheaves $(E,\phi^{E})$ such that $E^{\vee\vee}\cong G$, the framing $\phi^{E}$ is just the composition of the injection in the double dual with $\phi^G$ (and hence is fixed if we fix $\mathcal{G}$ and $E$) and the quotient $A=\sfrac{G}{E}$ is supported at $x_j$, with length over $x_j$ equal to $m_j$. We can decompose $A=\bigoplus A_j$ as sum of skyscraper sheaves supported at $x_j$ of length $m_j$; each of this $A_j$ may be seen as an Artin quotient of the trivial sheaf $O_X^{\oplus r}$ supported only at $x_j$ and of length $m_j$, hence the possible choice of $A_j$ are parametrized by the scheme $Quot_{x_j}(r,m_j)$; these schemes have been studied in the appendix of \cite{art:bara2}, we recall the following result:
\begin{prop} \label{prop:punctual_quot}
Let $X$ be a smooth projective surface, $x\in X$ a point and $r,d$ a pair of positive integers. Then the scheme $Quot_x(r,d)$ parametrizing Artin quotients of the trivial sheaf $O_X^{\oplus r}$ supported at $x$ of length $d$ is irreducible of dimension $rd-1$.
\end{prop}

We get as a corollary the following: 
\begin{cor}	\label{cor:fiber_hilbert-chow}
Let $\mathfrak{p}=(\mathcal{G},\sum_{j=1}^l m_j[x_j])$ be a point of $\mathcal{M}^{reg}_{X,D}(r,c,s)\times S^{n-s}X$.  Then $\pi_r^{-1}(\mathfrak{p})$ is isomorphic to the scheme $\prod_j Quot_{x_j}(r,m_j)$ as an irreducible closed subscheme of $\mathcal{M}_{X,D}(r,c,n)$ of dimension $\sum_j (rm_j-1)= r(n-s)-l$.
\end{cor}
	
\section{Nakajima operators on moduli spaces of framed sheaves}

From now on we fix the rank $r$ and the first Chern class $c$ and we consider the moduli spaces $\mathcal{M}_{X, D}(r,c,n)$ \emph{all together} by varying the second Chern class $n$. To emphasize this fact in the following we denote by $\mathcal{M}_{X,D}^{r,c}(n)$ the moduli space $\mathcal{M}_{X, D}(r,c,n)$ as if we consider $n$ a variable and $r, c$ constants.

Let $\mathrm{H}_{\bullet}(\mathcal{M}^{r,c}_{X,D}(n))$ be the homology group of $\mathcal{M}_{X,D}^{r,c}(n)$ with complex coefficients. Denote by $\mathbb{H}(X)$ the direct sum $\bigoplus_n \mathrm{H}_{\bullet}(\mathcal{M}^{r,c}_{X,D}(n))$ of all the homology groups of all the moduli spaces $\mathcal{M}^{r,c}_{X,D}(n).$

In this section for any $i\in \mathbb{Z}$ and $\alpha\in H_{\bullet}(X)$ we define operators 
\begin{equation*}
P_{\alpha}[i]: \mathbb{H}(X)\rightarrow \mathbb{H}(X)
\end{equation*}
that we shall call \emph{Nakajima operators} and in the following section we prove that they satisfy the Heisenberg commutator relations: in this way we define an action of a Heisenberg algebra associated to $\mathrm{H}_{\bullet}(X)$ on the space $\mathbb{H}(X)$ in a geometric way.

First for any $i,n\in \mathbb{Z}$ and $\alpha\in H_{\bullet}(X)$ we are going to construct operators 
\begin{equation*}
P_{\alpha}^n[i]:\mathrm{H}_{\bullet}(\mathcal{M}^{r,c}_{X, D}(n))\longrightarrow \mathrm{H}_{\bullet}(\mathcal{M}_{X, D}^{r,c}(n-i))
\end{equation*}
as \emph{correspondence} for a suitable subvariety of $\mathcal{M}^{r,c}_{X, D}(n-i)\times\mathcal{M}^{r,c}_{X, D}(n)$. Moreover, for technical reasons due to the non-projectivity of the moduli spaces $\mathcal{M}^{r,c}_{X, D}(n)$ we have to use the Borel-Moore homology groups of our moduli spaces (for technical details about Borel-Moore homology we refer to Appendix B in \cite{book:fulton-Young}).

Now we recall briefly the construction of correspondences. Let $M_1,M_2$ be two oriented manifolds of dimension $d_1, d_2$ respectively and $Z\subseteq M_1\times M_2$ an oriented submanifold such that the projection $p_2:Z\rightarrow M_2$ is a proper map. We denote by $\mathrm{H}^{lf}_{j}(M_a)$ the Borel-Moore homology group of $M_a$ for $a=1,2.$

We can define operators
\begin{eqnarray*}
\widetilde{\gamma}_{Z}:  \mathrm{H}^{lf}_j(M_1)&\longrightarrow &\mathrm{H}^{lf}_{j+\dim(Z)-d_1}(M_2)\\
	\eta & \longmapsto & p_{2*}(p_1^*\eta\cap [Z])
\end{eqnarray*}
and for any $\alpha\in \mathrm{H}^{lf}_{\bullet}(Z)$
\begin{eqnarray*}
\widetilde{\gamma}_\alpha: \mathrm{H}^{lf}_j(M_1) &\longrightarrow &\mathrm{H}^{lf}_{j+\deg(\alpha)-d_1}(M_2)\\
\eta & \longmapsto & p_{2*}(p_1^*\eta\cap \alpha)
\end{eqnarray*}
where the $\cap$-product is the one on Borel-Moore homology.

Moreover, using the fact that the singular homology group of $M_a$ is the direct limit of the Borel-Moore homology group of $M_a$ with respect to the compact subsets of $M_a$ (for $a=1, 2$), one defines operators between the singular homology groups (that we call \emph{correspondence})
\begin{equation*}
\gamma_{Z}:  \mathrm{H}_j(M_2)\longrightarrow  \mathrm{H}_{j+\dim(Z)-d_2}(M_1)
\end{equation*}
and for any $\alpha\in \mathrm{H}^{lf}_\bullet(Z)$
\begin{equation*}
\gamma_{\alpha}: \mathrm{H}_j(M_2)\longrightarrow \mathrm{H}_{j+\deg(\alpha)-d_2}(M_1).
\end{equation*}

Now we want to study the composition of two such homomorphism: if $M_1,M_2,M_3$ are three oriented manifolds and $Z\subseteq M_1\times M_2$, $W\subseteq M_2\times M_3$ two oriented submanifolds such that the projections $Z\rightarrow M_2$, $W\rightarrow M_3$ are proper maps, we have the composition of the correspondences $\gamma_{Z}\circ\gamma_{W}: \mathrm{H}_\bullet(M_3)\rightarrow \mathrm{H}_\bullet(M_1)$: these turn out to be equal to the correspondence defined by the subvariety $A\subseteq M_1\times M_3$,
\begin{equation*}
A=p_{13}(p_{12}^{-1}(Z)\cap p_{23}^{-1}(W)),
\end{equation*}
where $p_{ab}$ are the projections $M_1\times M_2\times M_3 \rightarrow M_a\times M_b$.

If we have $\alpha\in \mathrm{H}_\bullet(Z)$ and $\beta\in \mathrm{H}_\bullet(W)$, the composition of the correspondences $\gamma_{\beta}$ and $\gamma_{\alpha}$ is given by the correspondence induced by the class $\epsilon\in \mathrm{H}_\bullet(A)$ defined by 
\begin{equation*}
\epsilon=p_{13*}(p_{12}^*\alpha\cap p_{23}^*\beta).
\end{equation*}

Now consider the moduli spaces $\mathcal{M}^{r,c}_{X, D}(n)$; let $n, i \in \mathbb{Z},i>0$ and consider the closed subvariety $P^n[i]\subseteq\mathcal{M}^{r,c}_{X, D}(n)\times\mathcal{M}^{r,c}_{X, D}(n+i)\times X$ defined by
\begin{equation*}
P^n[i]=\left\{(\mathcal{E}_1,\mathcal{E}_2,x)\in\mathcal{M}^{r,c}_{X, D}(n)\times\mathcal{M}^{r,c}_{X, D}(n+i)\times X \left\vert\begin{array}{ll}
	\bullet & \mathcal{E}_1^{\vee\vee}\cong\mathcal{E}_2^{\vee\vee},\\
	\bullet & \mathcal{E}_2\subseteq \mathcal{E}_1,\\
	\bullet & \mathrm{supp}\left(\sfrac{E_1}{E_2} \right) = \{x\}.
	\end{array}\right.\right\}
\end{equation*}

We want to define \emph{corrispondences} using the fundalmental class of $P^n[i]$. Through direct computations one can prove that the dimension of the closed subvariety $P^n[i]$ is $2rn+b+r+1$.

Note that if $\pi:\mathcal{M}^{r,c}_{X, D}(n)\times\mathcal{M}^{r,c}_{X, D}(n+i)\times X\rightarrow\mathcal{M}^{r,c}_{X, D}(n)\times\mathcal{M}^{r,c}_{X, D}(n+i)$ is the projection, then the restriction $\pi\vert_{P^n[i]}$ is a proper injective map and therefore $\pi(P^n[i])$ is a closed subvariety. 

Thus for any $\alpha\in \mathrm{H}_\bullet(X)$ we have a Borel-Moore homology class $P^n_\alpha[i]\in \mathrm{H}^{lf}_\bullet (\pi(P^n[i]))$ defined by 
\begin{equation*}
P^n_\alpha[i]=\pi_*(\tau^*\alpha\cap [P^n[i]])
\end{equation*}
where $\tau:\mathcal{M}^{r,c}_{X, D}(n)\times\mathcal{M}^{r,c}_{X, D}(n+i)\times X\rightarrow X$ is the projection on the surface. (With the previous notation: $P^n_\alpha[i]=\widetilde{\gamma}_{P^n[i]}(\alpha)$).

Since the restriction to $\pi(P^n[i])$ of the projection from $\mathcal{M}^{r,c}_{X, D}(n)\times\mathcal{M}^{r,c}_{X, D}(n+i)$ to $\mathcal{M}^{r,c}_{X, D}(n+i)$ is a proper map, the class $P_\alpha^n[i]$ induces a correspondence operator
\begin{equation*}
P^n_\alpha[i]:\mathrm{H}_j(\mathcal{M}^{r,c}_{X, D}(n+i))\longrightarrow \mathrm{H}_{j+\deg(\alpha)-2ri-2}(\mathcal{M}^{r,c}_{X, D}(n))
\end{equation*}
(With the previous notation this would be $\gamma_{P^n_\alpha[i]}$).

In a similar way we can define for $i\in \mathbb{Z}, i>0$ closed subvarieties
\begin{equation*}
P^n[-i]\subseteq \mathcal{M}^{r,c}_{X, D}(n)\times\mathcal{M}^{r,c}_{X, D}(n-i)
\end{equation*}
as before exchanging the roles of $\mathcal{E}_1$ and $\mathcal{E}_2$; these subvarieties have dimension $2rn+b-ri+1.$ We can repeat the previous construction and for any $\beta\in \mathrm{H}_\bullet(X)$ we get the operator
\begin{equation*}
P_\beta^n[-i]: H_j(\mathcal{M}^{r,c}_{X, D}(n-i))\longrightarrow H_{j+\deg(\beta)+2ri-2}(\mathcal{M}^{r,c}_{X, D}(n)).
\end{equation*}
For $i=0$ we set $P^n_\alpha[0]$ equal to the identity map for any $n\in \mathbb{Z}$ and $\alpha\in \mathrm{H}_\bullet(X)$. 

Now we can extend these operators to $\mathbb{H}(X)$: for any $i\in \mathbb{Z}$ and $\alpha\in \mathrm{H}_{\bullet}(X)$ define operators
\begin{equation*}
P_\alpha[i]:\mathbb{H}_{\bullet}(X)\longrightarrow \mathbb{H}_{\bullet}(X)
\end{equation*}
acting on every summand $\mathrm{H}_\bullet(\mathcal{M}^{r,c}_{X, D}(n))$ with $P^n_\alpha[i]$.

\section{Main Theorem}

Now we want to study the commutation relations between the operators $P_\alpha[i]$ we constructed. In particular we want to prove the following result:

\begin{thm}	\label{thm:main}
For any $i,j\in\mathbb{Z}, \alpha,\beta\in \mathrm{H}_\bullet(X)$ the following holds:
\begin{equation*}
	[P_\alpha[i],P_\beta[j]] = (-1)^{ri-1}ri \delta_{i+j,0} <\alpha,\beta> \mathrm{id}
\end{equation*}
\end{thm}
This provides a representation of the Heisenberg algebra generated by $\sfrac{P_\alpha[i]}{(-1)^{ri-1}r}$ and $P_\beta[j]$ on $\mathbb{H}(X).$

The first step to get this result is the following:
\begin{prop}	\label{prop:quasi_main}
There exist constants $c_{r,i,n}$ such that for any $v\in \mathrm{H}_\bullet(\mathcal{M}^{r,c}_{X, D}(n))$ we have:
\begin{equation*}
[P_\alpha[i],P_\beta[j]] v= c_{r,i,n} \delta_{i+j,0} <\alpha,\beta> v 
\end{equation*}
\end{prop}
A proof of this proposition may be given directly by investigating the cycles $P^{n-j}_{\alpha}[i]P^{n}_{\beta}[j]$ and $P^{n-i}_{\beta}[j]P^{n}_{\alpha}[i]$ in $\mathcal{M}^{r,c}_{X, D}(n-i-j)\times\mathcal{M}^{r,c}_{X, D}(n)$, and observing that they act in the same way whenever $i\neq j$, while when $i=j$ they differ by a constant times the intersection product of the cycles; in \cite{book:naka} one can find the details of this proof for the case when $r=1$, i.e. the Hilbert scheme of points, and this directly generalize to our case.

\begin{lem}
The constants $c_{r,i,n}$ of the previous proposition are independent of $n$. (And so, from now on we shall denote them simply by $c_{r,i}$)
\end{lem}
\proof
Let $i,k\in\mathbb{Z}$ be such that $k\neq \pm i$, $\alpha,\beta,\gamma\in \mathrm{H}_\bullet(X)$ and  $v\in \mathrm{H}_\bullet(\mathcal{M}^{r,c}_{X, D}(n))$. 
Since $P_{\gamma}[k]v\in H_\bullet(\mathcal{M}^{r,c}_{X, D}(n-k))$, we get
\begin{eqnarray*}
[P_\gamma[k],[P_\alpha[i],P_{\beta}[-i]]]v &=& P_\gamma[k][P_\alpha[i],P_{\beta}[-i]]v - [P_\alpha[i],P_{\beta}[-i]] P_\gamma[k]v=\\
&=& <\alpha,\beta> (c_{r,i,n}-c_{r,i,n-k})v 
\end{eqnarray*}
But $k\neq \pm i$, hence $[P_\gamma[k],[P_\alpha[i],P_\beta[j]]]=0$, and we can choose $\alpha,\beta,\gamma,k,v$ so that $<\alpha,\beta>\neq 0$ and $P_\gamma[k]v\neq 0.$
	
So $c_{r,i,n}=c_{r,i,n+k}$; for $i\neq 0$, we can take $k=1$ and have $c_{r,i,n}=c_{r,i,n+1}$ that implies the thesis; for $i=1$ taking $k=2,3$ we have that $c_{r,1,n}=c_{r,1,n+2}=c_{r,1,n+3}$, and again we have the thesis.
\endproof
	
\section{Calculation of constants}

In this section we compute explicitly the constants $c_{r,i}$: following Grojnowski and Nakajima, we find an isomorphism of the algebra of symmetric functions with some subspaces of $\mathbb{H}_\bullet(X)$, such that the Newton polynomials correspond to the Nakajima operators $P_{[C]}[i]$ for some smooth connected curves $C\subseteq X$, $C\neq D.$ In the following section we recall some facts about the algebra of symmetric functions (for more informations see chapter one in \cite{book:mac}).

\textbf{Notations:} We use two different conventions for $\mu$ partition of a natural number $k$: if we write $\mu=(\mu_1,\ldots,\mu_t)$ we mean that $\sum_a \mu_a =k$ and assume that $\mu_a\geq\mu_{a+1}> 0$; moreover we call \emph{length} $\vert \mu \vert$ of the partition $\mu$ the natural number $t.$ If we write $\mu=(1^{m_1},\ldots,k^{m_k})$ we mean that $m_l=\#\{a\,\vert\,\mu_a=l\}$.

\subsection*{Symmetric Functions} 

We call \emph{algebra of symmetric polynomial in} $N$ \emph{variables} the subspace $\Lambda_N$ of $\mathbb{C}[x_1,\ldots,x_N]$ stable by the action of the $N$-th group of permutations $\sigma_N$. We have that $\Lambda_N$ is a graded ring:
\begin{equation*}
\Lambda_N=\bigoplus_n \Lambda_N^n
\end{equation*}
where $\Lambda_N^n$ is the ring of homogeneous symmetric polynomials in $N$ variables of degree $n$ (together with the zero polynomial).

For any $M>N$ we have morphisms $\rho_{MN}:\Lambda_M\rightarrow \Lambda_N$ mapping the variables $x_{N+1},\ldots,x_M$ to zero. Moreover the morphisms $\rho_{MN}$ preserve the grading, hence we can define $\rho_{MN}^n:\Lambda^n_M\rightarrow \Lambda^n_N$; this allows us to define
\begin{equation*}
\Lambda^n:= \lim_{\stackrel{\leftarrow}{N}} \Lambda^n_N
\end{equation*}
and the \emph{algebra of symmetric functions in infinitely many variables} as
\begin{equation*}
\Lambda:=\bigoplus_n \Lambda^n.
\end{equation*}

Now we want to define a basis for $\Lambda$, to do this we start by defining a basis in $\Lambda_N$. Let $\mu=(\mu_1,\ldots,\mu_t)$ be a partition with $t\leq N$, we define the polynomial
\begin{equation*}
m_\mu(x_1,\ldots,x_N) =\sum_{\tau\in \sigma_N(\mu)} x_1^{\tau_1}\cdots x_N^{\tau_N}
\end{equation*}
where we set $\mu_j=0$ for $j=t+1,\ldots, N$. The polynomial $m_\mu$ is symmetric, moreover the set of all $m_\mu$ for all the partitions $\mu$ with $\vert \mu \vert\leq N$ is a basis of $\Lambda_N$. If we consider the set of all $m_\mu$ for all the partitions $\mu$ with $\vert \mu \vert\leq N$ and $\sum_i \mu_i=n$, then this set is a basis of $\Lambda^n_N$.

Since for $M>N\geq t$ we have $\rho^n_{MN}(m_\mu(x_1,\ldots,x_M))=m_\mu(x_1,\ldots,x_N)$, by using the definition of inverse limit we can define the \emph{monomial symmetric functions} $m_\mu.$ By varying $\mu$ partition of $n$, these polynomials form a basis for $\Lambda^n.$

Now we want do define particular families of symmetric functions. Let $n \in \mathbb{N}, n \geq 1$,  we define the \emph{elementary symmetric function} $e_n$ as 
\begin{equation*}
e_n= m_{(1^n)} = \sum_{i_1<\ldots<i_n} x_{i_1}\dots x_{i_n}
\end{equation*}
and we put $e_0=1$. For $\mu=(\mu_1,\ldots,\mu_t)$ partition, we set $e_\mu=e_{\mu_1}e_{\mu_2}\dots e_{\mu_t}$; one can prove that the set $\{e_\mu\}_{\mu}$ is a basis of $\Lambda$.

The generating function for the sequence $\{e_n\}_{n\in\mathbb{N}}$ is 
\begin{equation*}
E(z)=\sum_{n=0}^\infty e_n z^n = \prod_{i=1}^{\infty} (1+x_i z).
\end{equation*}

Next, the \emph{$n$-th complete symmetric function} is $h_n=\sum_{\vert\nu\vert=n} m_\nu$. Since we have the formal power series
\begin{equation*}
\frac{1}{1-x_iz}= 1 + x_iz+x_i^2z^2+\ldots
\end{equation*}
the sequence $\{h_n\}_{n\in\mathbb{N}}$ has generating function
\begin{equation*}
H(z)=\sum_{n=0}^\infty h_n z^n = \prod_{i=1}^\infty \frac{1}{1-x_i z}.
\end{equation*}
In particular we have $E(z)H(-z)=1$.

Finally, the \emph{$n$-th Newton polynomial} $p_n$ is $p_n = m_{(n)}= \sum_{i=1}^\infty x_i^n$. Its generating function is
\begin{equation*}
P(z)=\sum_{n=1}^\infty p_n z^{n-1} =\sum_{n,i=1}^\infty x_i^n z^{n-1}
	= \sum_{i=1}^\infty \frac{d}{dz} \log \left( \frac{1}{1-x_iz} \right) = \frac{d}{dz} \log H(z).
\end{equation*}

Hence we have the formulas:
\begin{eqnarray} \nonumber
	H(z) &=& \exp \left(\int P(z)\right) = \exp \left( \sum_{n=0}^\infty \frac{p_nz^n}{n} \right),\\ \label{eqn:gen_symm_function}
	E(z) &=& H(-z)^{-1} = \exp \left( \sum_{n=0}^\infty \frac{p_nz^n}{(-1)^{n-1} n} \right).
\end{eqnarray}

To conclude this section, we state the following lemma:
\begin{lem} \label{lem:repr_symm_function}
For any $i\in\mathbb{N}$ we have
\begin{equation*}
p_i m_\mu = \sum_{\nu} a_{\mu\nu} m_\nu
\end{equation*}
where the summation is over partitions $\nu$ of $\vert \mu\vert+i$ which are obtained as follows: add $i$ to a term in $\mu$,say $\mu_j$ (possibly zero),and then arrange it in descending order. The coefficient $a_{\mu\nu}$ is $\#\{l\:|\:\exists k \ \text{s.t.}\ \nu_l = \mu_k+i\} $.
\end{lem}

\subsection*{Computation of $c_{r,i}$}

We need to introduce some subvarieties of our moduli spaces: let $C$ be a smooth connected curve in $X$, $C\neq D$, $s,k\in \mathbb{N}$; set
\begin{equation*}
L_C^k(s)^\circ = \left\{\mathcal{F}=(F, \phi^F)\in\mathcal{M}^{r,c}_{X,D}(s+k)\ \left\vert \begin{array}{ll}
	\bullet & \mathcal{F}\subset \mathcal{H} \mbox{ for some } \mathcal{H}\in\mathcal{M}^{r,c}_{X,D}(s),\\
	 \bullet & \sfrac{H}{F} \mbox{ is an Artin sheaf of length } k,\\
	 \bullet & \mathrm{supp}(\sfrac{H}{F})\subseteq C.
	 \end{array}\right.\right\}
\end{equation*}
and $L_C^k(s)$ its closure in $\mathcal{M}^{r,c}_{X,D}(s+k)$; for $\mu$ partition of $k$, set
\begin{equation*}
L_C^\mu(s)^\circ = \left\{\mathcal{F}=(F, \phi^F)\in\mathcal{M}^{r,c}_{X,D}(s+k)\ \left\vert \begin{array}{ll}
	 \bullet & \mathcal{F}\subset \mathcal{H} \mbox{ for some } \mathcal{H}\in\mathcal{M}^{r,c}_{X,D}(s),\\
	 \bullet & \mathrm{supp}(\sfrac{H}{F})=\{x_1,\ldots,x_t\}\subseteq C,\\ 
	 \bullet & \mathrm{length}_{x_a}(\sfrac{H}{F})=\mu_a  \mbox{ for } a=1, 2, \ldots, t.
	\end{array}\right.\right\}
\end{equation*}
and $L_C^\mu(s)$ its closure in $\mathcal{M}^{r,c}_{X,D}(s+k)$.

Now let $\mathcal{G}=(G,\phi^G)\in\mathcal{M}^{reg}_{X,D}(r,c,s)$ and set
\begin{equation*}
L_C^k(\mathcal{G})^\circ = \left\{\mathcal{F}=(F,\phi^F)\in\mathcal{M}^{r,c}_{X,D}(s+k)\ \left\vert \begin{array}{ll}
	 \bullet &  \mathcal{F}^{\vee\vee}= \mathcal{G},\\
	 \bullet & \mathrm{supp}(\sfrac{G}{F})\subseteq C.
	\end{array}\right.\right\}
\end{equation*}
and $L_C^k(\mathcal{G})$ its closure in $\mathcal{M}^{r,c}_{X,D}(s+k)$; for $\mu$ partition of $k$ set
\begin{equation*}
L_C^\mu(\mathcal{G})^\circ = \left\{\mathcal{F}=(F,\phi^F)\in\mathcal{M}^{r,c}_{X,D}(s+k)\ \left\vert \begin{array}{ll} 
	\bullet & \mathcal{F}^{\vee\vee}= \mathcal{G},\\
	\bullet & \mathrm{supp}(\sfrac{G}{F})=\{x_1,\ldots,x_t\}\subseteq C, \\ 
	\bullet & \mathrm{length}_{x_a}(\sfrac{G}{F})=\mu_a  \mbox{ for } a=1, 2, \ldots, t.
	\end{array}\right.\right\}
\end{equation*}
and $L_C^\mu(\mathcal{G})$ its closure. 

The first result we have is similar to that of proposition 4.5 in \cite{art:bara2}:
\begin{prop}\label{prop:irridu} 
The irreducible components of	$L_C^k(s)$ are the closed subvarieties $L_C^\mu(s)$ for $\mu$ a partition of $k$, each of which has dimension $2rs+rk$.
	
For any $\mathcal{G}\in\mathcal{M}^{reg}_{X,D}(r,c,s)$, the irreducible components of $L_C^k(\mathcal{G})$ are the closed subvarieties $L_C^\mu(\mathcal{G})$, for $\mu$ a partition of $k$, each of which has dimension $rk$.
\end{prop}

Now we remark the following: let $\alpha=[C]\in H_2(X)$ and consider the corresponding operator
\begin{equation*}
P_{[C]}[-i]: \mathrm{H}_j(\mathcal{M}^{r,c}_{X,D}(n))\rightarrow \mathrm{H}_{j+2ri}(\mathcal{M}^{r,c}_{X,D}(n+i))
\end{equation*}
we have that the subspace of $\mathbb{H}(X)$ of \emph{middle degrees}, i.e. $\bigoplus_{n} \mathrm{H}_{2rn}(\mathcal{M}^{r,c}_{X,D}(n))$, is preserved by the action of $P_{[C]}[-i]$.

Moreover, since $P_{[C]}[-i]$ `changes the sheaves only at a point of $C$ adding $i$ to the length of a zero-dimensional quotient', if $\mu$ is a partition of $k$ we have that $P_{[C]}[-i][L_C^\mu(\mathcal{G})]\in \mathrm{H}_{2r(k+i)}(L_C^{k+i}(\mathcal{G}))$; the next proposition calculates the coefficients of this element with respect to the basis of $\mathrm{H}_{2r(k+i)}(L_C^{k+i}(\mathcal{G}))$ consisting of the classes of the irreducible components of $L_C^{k+i}(\mathcal{G}).$

\begin{prop}\label{prop:repr_symm}
For any $s,k\in\mathbb{N}$, $\mu$ a partition of $k$, $C\subseteq X$ smooth connected curve, $C\neq D$, the following holds:
\begin{equation*}
P_{[C]}[-i] [L_C^\mu(s)] =\sum_{\nu \,\mathrm{partition}\, \mathrm{of}\, \vert\mu\vert+i} a_{\mu\nu} [L_C^\nu(s)]
\end{equation*}
where the constants $a_{\mu\nu}$ are the same as in lemma \ref{lem:repr_symm_function}.
	
For any $\mathcal{G}\in\mathcal{M}_{X,D}(r,c,s)$ we get
\begin{equation*}
P_{[C]}[-i] [L_C^\mu(\mathcal{G})]= \sum_{\nu \,\mathrm{partition}\, \mathrm{of}\, \vert\mu\vert+i} [L_C^\nu(\mathcal{G})]
\end{equation*}
with the same constants $a_{\mu\nu}$.
\end{prop}

The proof of this proposition is given in section \ref{sec:prop_pietro}.

Let $[vac]$ be the class of the point $\mathcal{G}$ and $[VAC]$ the fundamental class of $\mathcal{M}^{r,c}_{X,D}(s)$. Let $L(\mathcal{G})$ (resp. $L(s)$) be the subspaces of $\mathbb{H}(X)$ generated by the classes $[L^\mu_C(\mathcal{G})]$ and $[vac]$ (resp. $[L^\mu_C(s)]$ and $[VAC]$). We define a $\mathbb{C}$-linear isomorphism from $L(\mathcal{G})$ (resp. $L(s)$) onto the algebra $\Lambda$ of symmetric functions in infinitely many variables, sending $[vac]$ (resp. $[VAC]$) to $1\in \Lambda$ and $[L^\mu_C(\mathcal{G})]$ (resp. $[L^\mu_C(s)]$) to the monomial symmetric function $m_\mu$.

By proposition \ref{prop:repr_symm}, the Nakajima operator $P_{[C]}[-i]$ corresponds under this isomorphism to multiplication by the i-th Newton function $p_i\in \Lambda$.

By this isomorphism, the elementary symmetric function $e_n$ corresponds to $[L_C^{(1^n)}(\mathcal{G})]$ (resp. $[L_C^{(1^n)}(s)]$) and the equation \ref{eqn:gen_symm_function} gives the following equalities:
\begin{eqnarray} \label{eqn:repr_symm_1}
\sum_{n=0}^\infty z^n[L_C^{(1^n)}(s)] & = &  \exp\left( \sum_{i=1}^\infty \frac{z^i P_{[C]}[-i]}{(-1)^{i-1} i} \right) [VAC],\\ \label{eqn:repr_symm_2}
\sum_{n=0}^\infty z^n[L_C^{(1^n)}(\mathcal{G})] & = &\exp\left( \sum_{i=1}^\infty \frac{z^i P_{[C]}[-i]}{(-1)^{i-1} i} \right) [vac].
\end{eqnarray}

Now we need the following proposition, whose proof will be given in section \ref{sec:prop_franco}:

\begin{prop}\label{prop:intersection_fix_non_fix}
Let $L$ be a very ample divisor in $X$. Let $C$ and $C'$ be two smooth connected curves different from $D$, such that $C\cap C'\cap D=\emptyset.$ Assume that $C$ and $C'$ meets trasversally and are in $\vert L \vert$. Then	the following holds:
\begin{equation*}
	\sum_{n=0}^\infty z^{2n} <[L_C^{(1^n)}(\mathcal{G})],[L_{C'}^{(1^n)}(s)]>= (1-(-1)^r z^2 )^{r<[C],[C']>}.
\end{equation*}
\end{prop}

We can finally conclude the proof of theorem \ref{thm:main}, via the same computations of the proof of theorem 4.1 in \cite{art:bara2}:
\begin{prop}
The constants $c_{r,i}$ of proposition \ref{prop:quasi_main} satisfy
\begin{equation*}
	c_{r,i} = (-1)^{ri-1} ri.
\end{equation*}
\end{prop}
\proof
For any $C$ smooth connected curve in $X$, $C\neq D$, set 
\begin{equation*}
C_+(z) = \sum_{i=1}^\infty \frac{P_{[C]}[i]z^i}{(-1)^{i-1}i},\; C_-(z) = \sum_{i=1}^\infty \frac{P_{[C]}[-i]z^i}{(-1)^{i-1}i}. 
\end{equation*}

Since $P_{[C]}[i]$ is the adjoint to $P_{[C]}[-i]$ with respect to the intersection product on $\mathbb{H}(X)$, $C_+(z)$ is the adjoint to $C_-(z)$ (with respect to the intersection product extended by linearity to power series).

Now, by proposition \ref{prop:intersection_fix_non_fix} and formulas \ref{eqn:repr_symm_1} and \ref{eqn:repr_symm_2} we have the equalities
\begin{eqnarray*}
 (1-(-1)^r z^2)^{<[C],[C']>} & = & \sum_{n=0}^\infty z^{2n} <[L_C^{(1^n)}(\mathcal{G})],[L_{C'}^{(1^n)}(s)]> = \\
&=& <\sum_{n=0}^\infty z^n [L_C^{(1^n)}(\mathcal{G})],\sum_{n=0}^\infty z^n [L_{C'}^{(1^n)}(s)]> = \\
&=& <\exp (C_-(z))[vac],\exp(C'_-(z))[VAC]> = \\
&=& <\exp(C'_+(z))\exp(C_-(z))[vac],[VAC]> = \\
&=& <\exp(C_-(z)) \left[\exp (-\mathrm{ad} C_-(z)))\left(\exp(C'_+(z)) \right)\right][vac],[VAC]>
\end{eqnarray*}
where the last equality follows from the general equality between operators $\exp(-A)B\exp(A) = \exp(-\mathrm{ad} A)(B)$ (we put $A=C_-(z)$ and $B=\exp(C'_+(z)$).

An explicit computation shows that 
\begin{equation*}
[C_-(z),\exp(C'_+(z))] = - \Phi(z) \exp(C'_+(z)) 
\end{equation*}
with $\Phi(z)=\sum_{n=0}^\infty \frac{c_{r,n}}{n^2} <[C],[C']> z^{2n}$; from this and the fact that $\exp(-\mathrm{ad} A)(B) = 1 - \mathrm{ad} A (B) + \ldots$ one gets that
\begin{equation*}
\exp (-\mathrm{ad} C_-(z)))\left(\exp(Cd_+(z)) \right) = \exp(\Phi(z))\exp(Cd_+(z)),
\end{equation*}
and so 
\begin{eqnarray*}
(1-(-1)^r z^2)^{<[C],[C']>} &=&\exp(\Phi(z)) <\exp(C_-(z))\exp(C'_+(z)) [vac],[VAC]> = \\
&=& \exp(\Phi(z)) <\exp(C_-(z))[vac],[VAC]> = \\
&=& \exp(\Phi(z)).
\end{eqnarray*}
where the second of this equalities follows from the fact that since $i>0$, $P_{[C']}[i][vac]=0$, while the last follows from the fact that every $P_{[C]}[-i]$ maps $[vac]$ in the orthogonal complement of $[VAC]$.

So finally we get
\begin{equation*}
r<[C],[C']> \log (1-(-1)^r z^2) = \sum_{n=1}^\infty \frac{c_{r,n}}{n^2} <[C],[C']> z^{2n} 
\end{equation*}
from which, developing in power series the left hand side, we get $c_{r,n}=(-1)^{rn-1}rn$, that is the thesis.
\endproof
	
\section{Proof of proposition \ref{prop:repr_symm}}\label{sec:prop_pietro}

Now we want to give a proof of proposition \ref{prop:repr_symm}; we explain how to do with $L_C^\mu(\mathcal{G})$, the case of $L_C^\mu(s)$ is analogous.

First of all we need the following lemma:
\begin{lem}	\label{lem:generic}
Let $\lambda$ be a partition of $k.$ For $\mathcal{F}$ generic in $L_C^\lambda(\mathcal{G})$ (i.e. in a open dense subset), around each $x_j$ point of the support of $\sfrac{G}{F}$ exist coordinates $(\eta_j,\zeta_j)$	such that 
\begin{equation*}
\left(\frac{G}{F}\right)_{x_j} \cong \frac{\mathbb{C}[\zeta_j]}{(\zeta_j)^{\lambda_j}}.
\end{equation*}
\end{lem}
\proof
It follows from the first paragraph of section 5 in \cite{art:bara2}.
\endproof

As we have already seen, $P_{[C]}[-i]([L_C^\mu(\mathcal{G})])$ is an element of maximal degree in $H_\bullet(L_C^{k+i}(\mathcal{G}))$. Since $L_C^\nu(\mathcal{G})$ are the irreducible components of $L_C^{k+i}(\mathcal{G})$  for $\nu$ partition of $k+i$, $P_{[C]}[-i]([L_C^\mu(\mathcal{G})])$ decomposes as a sum $\sum_{\nu} b_{\mu\nu}[L_C^\nu(\mathcal{G})]$ for some coefficients $b_{\mu\nu}.$ 

Put $n=k+s.$ To compute the coefficients $b_{\mu\nu}$, consider the intersection
\begin{equation*}
P_{[C]}^n[-i]\cap p_n^{-1}(L_C^\mu(\mathcal{G}))\cap p_{n+i}^{-1}(L_C^\nu(\mathcal{G}))
\end{equation*}
where the projection morphisms are:
\begin{equation*}
\xymatrix{
&\mathcal{M}_{X,D}^{r,c}(n+i)\times\mathcal{M}_{X,D}^{r,c}(n)\ar_{p_{n+i}}[dl]\ar^{p_{n}}[dr] & \\
\mathcal{M}_{X,D}^{r,c}(n+i) & & \mathcal{M}_{X,D}^{r,c}(n) }
\end{equation*}

We prove the proposition:
\begin{prop}\label{prop:transversal_intersection}
Let $\mathcal{F}_2\in L_C^\mu(\mathcal{G}), \mathcal{F}_1\in L_C^\nu(\mathcal{G})$ be such that $(\mathcal{F}_1,\mathcal{F}_2)\in P_{[C]}^n[-i]$ is a smooth point and such that they are generic in the sense of lemma \ref{lem:generic}. Then the intersection of the tangent spaces
\begin{equation*}
W= T_{(\mathcal{F}_1,\mathcal{F}_2)} p_n^{-1}(L_C^\mu(\mathcal{G}))
		\cap T_{(\mathcal{F}_1,\mathcal{F}_2)} P_{[C]}^n[-i]
\end{equation*}
is isomorphic via $p_{n+i\ *}$ to the tangent space $T_{\mathcal{F}_1} L_C^\nu(\mathcal{G}).$
\end{prop}

This result implies that the coefficients $b_{\mu\nu}$ may be computed by counting the number of points in a generic fiber of the morphism
\begin{equation*}
p_{n+i}\vert_{P_{[C]}[-i]\cap p_n^{-1}(L_C^\mu(\mathcal{G}))\cap p_{n+i}^{-1}(L_C^\nu(\mathcal{G}))}.
\end{equation*}

Fix a generic (in the sense of lemma \ref{lem:generic}) $\mathcal{E}_2\in L_C^\nu(\mathcal{G})$ such that the quotient sheaf $A_2=\sfrac{G}{E_2}$ is of the form stated in lemma \ref{lem:generic}. For a generic framed sheaf $\mathcal{E}_1\in L_C^\mu(\mathcal{G})$ such that $(\mathcal{E}_1,\mathcal{E}_2)\in P_{[C]}[-i]$ the quotient $A_1=\sfrac{G}{E_1}$ satisfy the equality $(A_1)_{x_i}=\frac{\mathbb{C}[\zeta_i]}{(\zeta_i)^{\mu_i}}$, with $\mu_i=\nu_i$ for all $i$ except for exactly one $i_0$, for which $\mu_{i_0}=\nu_{i_0}-i.$ So the number of generic $\mathcal{E}_1\in L_C^\mu(\mathcal{G})$ such that $(\mathcal{E}_1,\mathcal{E}_2)\in P_{[C]}[-i]$ is exactly $a_{\mu\nu}$, thus $b_{\mu\nu}=a_{\mu\nu}.$

To prove proposition \ref{prop:transversal_intersection} we need to give a description of the tangent spaces 
to $L_C^\mu(\mathcal{G})$ and $P_{[C]}[-i]$.

\subsection*{Tangent to $L_C^\lambda(\mathcal{G})$}

Let $\lambda$ be a partition of $k.$ Fix $\mathcal{F}=(F,\phi^F)\in L_C^\lambda(\mathcal{G}).$ Remark that since $\mathcal{F}^{\vee\vee}=\mathcal{G}$, the framing $\phi^F$ is just the composition of the inclusion $F\hookrightarrow G$ with the framing $\phi^G$ of $\mathcal{G}.$ So the deformations of $\mathcal{F}$ in $L_C^\lambda(\mathcal{G})$ (that in general would be deformations of $F$ and $\phi^F$) are uniquely determined by deformations of the sheaf $F$ with respect to the functor $\mathscr{F}(\underline{L_C^\lambda(\mathcal{G})})$, where $\mathscr{F}$ is the `framing-forgetting' morphism of functors $\underline{\mathcal{M}_{X,D}^{r,c}}(n) \rightarrow \mathscr{S}$, where $\mathscr{S}$ is the controvariant functor from the category of schemes to the category of sets that associates to every scheme $S$ the set consisting of isomorphism classes of flat families of torsion free sheaves over $X$ parametrized by $S.$

The deformations of the sheaf $F$ have been studied in \cite{art:bara2}, here we recall that study. 

First note that the functor $\mathscr{F}(\underline{L_C^\lambda(\mathcal{G})})$ is a subfunctor of the Quot-scheme functor $\underline{Quot}_X(G,\vert\lambda\vert).$ It is well known that deformations at a point $F$ in $Quot_X(G,\vert\lambda\vert)$ are parametrized by $\mathrm{Hom}(F,A)$, where $A=\sfrac{G}{F}$, that can be seen as a subspace of $\mathrm{Ext}^1(F,F)$ via the connection morphism that comes applying the left-exact functor $\mathrm{Hom}(F,\bullet)$ to the short exact sequence $0\rightarrow F\rightarrow G\rightarrow A\rightarrow 0.$

Now, for $\mathcal{F}=(F,\phi^F)\in L_C^\lambda(\mathcal{G})$ we have the decomposition $A=\bigoplus_j A_j$ with $A_j$ a skyscraper sheaf, supported at $x_j\in C$ of length $\lambda_j$: this gives a decomposition $\mathrm{Hom}(F,A)=\bigoplus_j\mathrm{Hom}(F,A_j)$. If we assume that $\mathcal{F}$ is generic, around each $x_j$ we have a neighborhood in which we have an inclusion $G(-\lambda_j C)\hookrightarrow F$ such that the diagram
\begin{equation} \label{diag:support_m}
		\xymatrix{
			G(-\lambda_j C)\ar@{^{(}->}[r]\ar@{=}[d] & F\ar@{^{(}->}[d]\\
			G(-\lambda_j C)\ar@{^{(}->}[r] & G
		}
\end{equation}
commutes; it turns out that deformations of $F$ with respect to $\mathscr{F}(\underline{L_C^\lambda(\mathcal{G})})$ are the one that preserves these diagrams. Now with deformation theory and diagram chasing one can prove
\begin{prop}
	The set of deformations of a sheaf $F$ with respect to the functor $\mathscr{F}(\underline{L_C^\lambda(\mathcal{G})})$ that preserve the diagram \ref{diag:support_m} are parametrized by
\begin{equation*}
\bigoplus_{j=1}^t \mathrm{Hom}\left(\frac{F}{G(-\lambda_j C)},
		\left(\frac{G}{F}\right)_{x_j}\right).
\end{equation*}
\end{prop}

So the tangent space of $L_C^\lambda(\mathcal{G})$ at a generic framed sheaf $\mathcal{F}$ is 
\begin{equation*}
T_{\mathcal{F}}(L_C^\lambda(\mathcal{G}))=\bigoplus_{j=1}^t \mathrm{Hom}\left(\frac{F}{G(-\lambda_j C)},
		\left(\frac{G}{F}\right)_{x_j}\right).
\end{equation*}

Remark that this is a subspace of $\mathrm{Ext}^1(F,F(-D))$ via
\begin{equation*}
\bigoplus_{j=1}^t \mathrm{Hom}\left(\frac{F}{G(-\lambda_j C)},
		\left(\frac{G}{F}\right)_{x_j}\right)\hookrightarrow \mathrm{Hom}(E,A)=\mathrm{Hom}(E,A(-D))\hookrightarrow \mathrm{Ext}^1(F,F(-D)).
\end{equation*}

\subsection*{Tangent to $P_{[C]}^n[-i]$}

We recall that $P_{[C]}^n[-i]$ is the subset of $\mathcal{M}_{X,D}^{r,c}(n+i)\times\mathcal{M}_{X,D}^{r,c}(n)$ consisting of pairs $(\mathcal{F}_1,\mathcal{F}_2)$ such that
\begin{itemize}
	\item[(a)] $\mathcal{F}_1^{\vee\vee}\cong \mathcal{F}_2^{\vee\vee}$;
	\item[(b)] $F_1\subseteq F_2$;
	\item[(c)] $\mathrm{supp}(\sfrac{F_2}{F_1})=\{x\}$ with $x\in C$ and $\mathrm{length}(\sfrac{F_2}{F_1})=i.$
\end{itemize}
Condition (a) implies that the framing $\phi^{F_1}$ of $F_1$ is uniquely determined by the inclusion morphism and by the framing on $F_2.$ So the deformations of $(\mathcal{F}_1,\mathcal{F}_2)$ are one to one with triples $(F_1', F_2', \phi')$ where $F_a'$ is a deformation of $F_a$ for $a=1,2$ and $\phi'$ is a deformation of the framing $\phi^{F_2}.$ Moreover, note both pairs $(F_1',\phi')$ and $(F_2',\phi')$ defines deformations of $\mathcal{F}_1$ and $\mathcal{F}_2$, hence define elements of $\mathrm{Ext}^1(F_1, F_1(-D))$ and $\mathrm{Ext}^1(F_2, F_2(-D))$ respectively.

The pair $(F_1',F_2')$ is a deformation of $(F_1,F_2)$  with respect to the functor $\mathscr{F}(\underline{P_{[C]}^n[-i]}).$ This kind of deformations are studied by Baranovsky in \cite{art:bara2}, here we recall his study.

Conditions (b) and (c) are equivalent to having in a neighborhood of $x$ the commutative diagram
\begin{equation*} 
		\xymatrix{
			F_2(-i C)\ar@{^{(}->}[r]\ar@{=}[d] & F_1\ar@{^{(}->}[d]\\
			F_2(-i C)\ar@{^{(}->}[r] & F_2
		}
\end{equation*}
where the lower arrow is the multiplication by $i$ times a generator of the curve $C.$ So a deformation of $(F_1,F_2)$ with respect to $\mathscr{F}(\underline{P_{[C]}^n[-i]})$ is a pair of sheaves $(F_1',F_2')$ over $X\times \mathrm{Spec}\left(\sfrac{\mathbb{C}[\varepsilon]}{\varepsilon^2}\right)$, flat over $\mathrm{Spec}\left(\sfrac{\mathbb{C}[\varepsilon]}{\varepsilon^2}\right)$, with fibers over the closed point isomorphic respectively to $F_1$ and $F_2$ satisfying the following conditions:
\begin{itemize}
 \item[(i)] there is an inclusion $F_1'\hookrightarrow F_2'$;
 \item[(ii)] in a neighborhood of $x$ we have an inclusion $F_2'(-i C)\hookrightarrow F_1'$;
 \item[(iii)] in this same neighborhood we have a commutative diagram 
	\begin{equation}\label{eq:condiii} 
		\xymatrix{
			F_2'(-i C)\ar@{^{(}->}[r]\ar@{=}[d] & F_1'\ar@{^{(}->}[d]\\
			F_2'(-i C)\ar@{^{(}->}[r] & F_2'
		}
	\end{equation}
	where the lower arrow is the multiplication by $i$ times a generator of the curve $C.$
\end{itemize}

Recall that for $a=1,2$, the sheaves $F_a'$ fit in an exact sequence
\begin{equation*}
0\rightarrow F_a\rightarrow F_a'\rightarrow F_a\rightarrow 0
\end{equation*}
that corresponds to an element $v_a\in\mathrm{Ext}^1(F_a,F_a).$ By diagram chasing, we get that condition (i) is equivalent to the fact that the images of $v_1$ and $v_2$ in $\mathrm{Ext}^1(F_1,F_2)$ coincide.

Condition (ii) may be treated in the same way: first of all notice that $v_2\in\mathrm{Ext}^1(F_2,F_2)$ defines an element $\tilde{v_2}\in\mathrm{Ext}^1(F_2(-i C),F_2(-i C))$. One can prove that condition (ii) is equivalent to the fact that the images of $v_1$ and $\tilde{v_2}$ in $\mathrm{Ext}^1(F_2(-i C),F_1)$ coincide.

Finally we come to condition (iii). Assume that $(F_1', F_2')$ satisfies conditions (i) and (ii) and label the morphisms:
\begin{equation*}
\tilde{\rho}:F_2'(-i C)\rightarrow F_1',\;\; \sigma:F_2'(-i C)\rightarrow F_2'
\end{equation*}
and $\rho$ the composition of $\tilde{\rho}$ with the inclusion $F_1'\rightarrow F_2'.$ This allows us to define a morphism $\tau:F_2(-i C)\rightarrow F_1$ that keeps track of the commutativity of the diagram \ref{eq:condiii}, i.e. such that $\tau$ vanish if and only if condition (iii) holds.

These three facts completely describe the deformations of pairs of sheaves $(F_1,F_2).$ So the tangent of $P_{[C]}^n[-i]$ at point $(F_1, F_2)$ is parametrized by pairs 
\begin{equation*}
(v_1,v_2)\in \mathrm{Ext}^1(F_1,F_1)\oplus \mathrm{Ext}^1(F_2,F_2)
\end{equation*}
such that the corresponding extensions satisfy conditions (i), (ii) and (iii) and both $v_a\in \mathrm{Ext}^1(F_a, F_a(-D))$ for $a=1,2.$

\subsection*{Proof of proposition \ref{prop:transversal_intersection}}

Let $w\in W.$ By previous descriptions $w=(v_1,v_2)\in  \mathrm{Ext}^1(F_1,F_1(-D))\oplus \mathrm{Ext}^1(F_2,F_2(-D))$ is such that the corresponding extensions satisfy conditions (i), (ii) and (iii) and moreover $v_2\in \bigoplus_{j=1}^t \mathrm{Hom}\left(\frac{F_2}{G(-\mu_j C)},(A_2)_{x_j}\right)$, where $A_2=\sfrac{G}{F_2}.$ We want to prove that $v_1\in T_{\mathcal{F}_1} L_C^\nu(\mathcal{G})$ and that $v_1=0$ only if $v_2=0.$ 

Since $\mathcal{F}_1^{\vee\vee}\cong \mathcal{F}_2^{\vee\vee}\cong \mathcal{G}$, we have that $v_1$ is an element of $\mathrm{Hom}(F_1,A_1)$, where $A_1$ is the quotient sheaf $\sfrac{G}{F_1}.$ We need to prove that actually $v_1\in\bigoplus_{j=1}^t \mathrm{Hom}\left(\frac{F_1}{G(-\nu_j C)},(A_1)_{x_j}\right).$

Using the inclusion morphism from $F_a$ to its double dual $G$, for $a=1,2$, and the inclusion morphism $F_1\rightarrow F_2$, we get the commutative diagram
\begin{equation}\label{diag:trasvers}
\xymatrix{
	\mathrm{Hom}(F_1,A_1) \ar[r]\ar[d] & 
		\mathrm{Ext}^1(F_1,F_1) \ar[d]\\
	\mathrm{Hom}(F_1,A_2) \ar[r] & 
		\mathrm{Ext}^1(F_1,F_2)\\
	\mathrm{Hom}(F_2,A_2) \ar[u] \ar[r]& 
		\mathrm{Ext}^1(F_2,F_2) \ar[u]\\
}
\end{equation}

\begin{lem}\label{lem:arrow}
The horizontal arrows of the previous diagram are injective morphisms.
\end{lem}
\proof
Consider the first row: it comes from the exact sequence
\begin{equation*}
0\rightarrow \mathrm{Hom}(F_1,F_1) \rightarrow 	\mathrm{Hom}(F_1,G) \rightarrow  \mathrm{Hom}(F_1,A_1) \rightarrow \mathrm{Ext}^1(F_1,F_1)
\end{equation*}
It suffices to show that the first arrow is an isomorphism, but it is trivial because by lemma \ref{lem:morph_framed} we have $\mathrm{Hom}(F_1,F_1) \cong \mathrm{End}(\mathbb{C}^r)\cong\mathrm{Hom}(F_1,G).$ By applying the same argument to the other arrows, we get the thesis.
\endproof

By condition (ii) the images of $v_1$ and $v_2$ in $\mathrm{Ext}^1(F_1,F_2)$ coincide. Since the diagram \ref{diag:trasvers} commutes and the second horizontal arrow is injective, the images of $v_1$ and $v_2$ actually coincide in $\mathrm{Hom}(F_1,A_2).$ Recall that the sheaves $A_{1}$ and $A_2$ are supported at points and differ only at $x$, so we can decompose $A_1$ and $A_2$ in the following way:
\begin{equation*}
	A_1=\bigoplus_{j} A_{x_j} \oplus B_1 \qquad 
	A_2 =\bigoplus_{j} A_{x_j} \oplus B_2
\end{equation*}
where $B_1,B_2$ are sheaves supported in $x$, and $\mathrm{length}(B_1)-\mathrm{length}(B_2)=i$ (remark that $B_2$ may be zero). So we can decompose the morphisms $v_1,v_2$ in
\begin{equation*}
v_1=\bigoplus_j v_{1,j} \oplus u_1\qquad v_2=\bigoplus_j v_{2,j} \oplus u_2
\end{equation*}
where $v_{a,j}: (F_a)_{x_j}\rightarrow A_{x_j}$ and $u_a:(F_a)_x\rightarrow B_a$ for $a=1,2.$ Since the images of $v_1,v_2$ in $\mathrm{Hom}(F_1,A_2)$ coincide and $F_1$ differs from $F_2$ only at $x$, we have that $v_{1,j}=v_{2,j}.$ Moreover at each point $x_j$ the morphism $v_{2,j}$ factors through $\sfrac{F_2}{G(-\mu_j C)}$, so $v_{1,j}$ factors through $\sfrac{F_1}{G(-\nu_j C)}$; so it remains to show only that $u_1$ factors through $\sfrac{F_1}{G(-\nu_{j_0} C)}$ knowing that $u_2$ factors through $\sfrac{F_2}{G(-\mu_{j_0} C)}$, and $\nu_{j_0}=\mu_{j_0}+i.$ Recall that in a neighborhood of $x$ we have an inclusion morphism $F_2(-i C)\hookrightarrow F_1$: this induces a diagram similar to \ref{diag:trasvers}:
\begin{equation*}
\xymatrix{
	\mathrm{Hom}(F_2(-i C),A_2(-i C)) \ar[r]\ar[d] & 
		\mathrm{Ext}^1(F_2(-i C),F_2(-i C)) \ar[d]\\
	\mathrm{Hom}(F_2(-i C),A_1) \ar[r] & 
		\mathrm{Ext}^1(F_2(-i C),F_1)\\
	\mathrm{Hom}(F_1,A_1) \ar[u] \ar[r]& 
		\mathrm{Ext}^1(F_1,F_1) \ar[u]
}\end{equation*}
with injective horizontal arrows. Condition (ii) assures that the images of $\tilde{v_2},v_1$ in $\mathrm{Hom}(F_2(-i C),A_1)$ coincide, hence it follows that $v_1$ factors through $\sfrac{F_1}{G(-\nu_{j_0} C)}$ with $\nu_{j_0}=\mu_{j_0}+i.$

Now assume that $v_1=0$, we want to show that this implies that $v_2=0$. From the commutativity of diagram \ref{eq:condiii} we get the commutativity of the diagram
\begin{equation*}\xymatrix{
	F_2(-i C)\ar[r]\ar^{\tilde{v_2}}[d] & F_1 \ar^{v_1}[d]\\
	A_2(-i C)\ar[r] & A_1
}\end{equation*}

The upper horizontal arrow in this diagram is injective by hypothesis and this implies that the lower one is injective too. Now from the fact that $v_1=0$ follows that $\tilde{v_2}=0$, hence $v_2=0$ and this concludes the proof of proposition \ref{prop:repr_symm}.

\section{Proof of proposition \ref{prop:intersection_fix_non_fix}}\label{sec:prop_franco}

In this section we provide a proof of proposition \ref{prop:intersection_fix_non_fix}; to achieve this, we compute the intersection number of the fundamental classes of the cycles $L_{C'}^{(1^n)}(s)$ and $L_C^{(1^n)}(\mathcal{G})$ in $\mathcal{M}^{r,c}_{X, D}(n+s)$. 

First note that any framed sheaf $\mathcal{F}$ in $L_C^{(1^n)}(\mathcal{G})$ corresponds to the kernel of a composition
\begin{equation*}
G\rightarrow G\vert_{C} \rightarrow A \rightarrow 0
\end{equation*}
where $A$ is an Artin sheaf of length $n$ supported on $C$. Thus $L_C^{(1^n)}(\mathcal{G})$ is isomorphic to $Quot_{C}^n(G\vert_{C})$, that is the Quot scheme that parametrizes the length $n$ Artin quotient sheaves of $G\vert_{C}.$

Let $\chi: G\vert_C\rightarrow A$ be an element in $Quot_{C}^n(G\vert_{C}).$ Denote by $F_1$ the kernel of $\chi$ (note that it is a locally free sheaf). By \textit{Serre's duality} we have
\begin{equation*}
\mathrm{Ext}^1(F_1, A)=\mathrm{Hom}(A\otimes F_1^{\vee}, \omega_C)=\mathrm{H}^{1}(C, A\otimes F_1^{\vee})
\end{equation*}
and $\mathrm{H}^{1}(C, A\otimes F_1^{\vee})=0$ because $A$ is a zero-dimensional sheaf, hence by proposition 2.2.8 of \cite{book:hl} we get that $Quot_{C}^n(G\vert_{C})$ is smooth at the point $A$. Thus $L_C^{(1^n)}(\mathcal{G})$ is smooth.

Moreover there exists an open subset $L_{C'}^{(1^n)}(s)^{\bullet}$ in $L_{C'}^{(1^n)}(s)$ isomorphic to a fiber bundle over $\mathcal{M}^{reg}_{X, D}(r,c,s)$ with fibers $L_{C'}^{(1^n)}(\mathcal{H})$ for $\mathcal{H}\in \mathcal{M}^{reg}_{X, D}(r,c,s).$

Finally we observe that $L^{(1^n)}_C(\mathcal{G})$ does not intersect $L_{C'}^{(1^n)}(s)\setminus L_{C'}^{(1^n)}(s)^\bullet$ because for any $\mathcal{F}\in L^{(1^n)}_C(\mathcal{G})$ the length of the Artin sheaf $\sfrac{F^{\vee\vee}}{F}$ is equal to $n$ while, on the other hand, if $\mathcal{F}\in L_{C'}^{(1^n)}(s)\setminus L_{C'}^{(1^n)}(s)^\bullet$, then $\mathrm{length}\left(\sfrac{F^{\vee\vee}}{F}\right)\geq n+1.$

Now we compute the intersection in the case where $X=\mathbb{C}\mathbb{P}^2$ and $C,C'$ are two distinct lines intersecting at $x\notin D$; at the end of the section we explain how the general case follows from this.

\subsection*{$X=\mathbb{C}\mathbb{P}^2$, $C,C'$ distinct lines}

Let $\mathcal{F}=(F, \phi^{F})\in L_C^{(1^n)}(\mathcal{G})\cap L_{C'}^{(1^n)}(s)$; then $\mathcal{F}^{\vee\vee} \cong \mathcal{G}$, the quotient $\sfrac{G}{F}$ is of length $n$ and supported on both $C$ and $C'$, hence it is supported only at $x$; moreover $\left(\sfrac{G}{F}\right)_x$ is a vector space of dimension $n$, and it is a quotient of $G_x$, wich is an $r$-dimensional vector space. So the intersection $L_C^{(1^n)}(\mathcal{G})\cap L_{C'}^{(1^n)}(s)$ is parametrized by the $n$-dimensional quotients of $G_x$, i.e. by the Grassmann variety $Gr(G_x,n)$ (in particular we observe that if $n>r$ the intersection is empty).
	
Now $L_C^{(1^n)}(\mathcal{G})$ and $L_{C'}^{(1^n)}(s)$	have complementary dimension in $\mathcal{M}^{r,c}_{X, D}(s+n)$, but their set-theoretic intersection has positive dimension; however we are dealing with this situation:
\begin{equation*}	
	\xymatrix{ 
	{Gr(G_x, n)} \ar[dr]^h \ar[d]_g \ar[r]^j & {L_{C'}^{(1^n)}(s)} \ar[d]^f \\
	{L_C^{(1^n)}(\mathcal{G})} \ar[r]_i      & {\mathcal{M}^{r,c}_{X, D}(n+s)}        }
\end{equation*}
where $i$ is a regular imbedding of codimension $2rs+rn+b$ and $j$ is a regular imbedding of codimension $2rs+n^2+b$ by theorem 17.12.1 in \cite{book:egaIV}. Thus we can apply the \textit{excess intersection formula} (see chapter 6 in \cite{book:fulton}) to get
\begin{equation}
i_{*}[L_C^{(1^n)}(\mathcal{G})] \cdot f_{*}[L_{C'}^{(1^n)}(s)] = h_{*} (c_{n(r-n)}(V) \cap [Gr(G_x,n)])
\end{equation}
where
\begin{equation*}
V=\frac{g^{*}N_{L_C^{(1^n)}(\mathcal{G})}(\mathcal{M}^{r,c}_{X, D}(n+s))}{N_{Gr(G_x,n)}(L_{C'}^{(1^n)}(s))}
\end{equation*}
is a locally free sheaf of rank $n(r-n)$. It is called \textit{excess bundle}.

The vector bundle $V$ arises also from the following exact sequence
\begin{equation*}
	\xymatrix{ 
	{0} \ar[r] & {T_{Gr(G_x,n)}} \ar[r] &  {g^{*}T_{L_C^{(1^n)}(\mathcal{G})}\oplus j^{*}T_{L_{C'}^{(1^n)}(s)}} \ar[r] & {h^{*}T_{\mathcal{M}^{r,c}_{X, D}(n+s)}} \ar[r] & {V} \ar[r] & {0}        }
\end{equation*}

Using this sequence, we obtain the following relation among the full Chern classes of the tangent bundles:
\begin{equation}\label{eq:chern}
c(V)=\frac{c(h^{*}T_{\mathcal{M}^{r,c}_{X, D}(n+s)})c(T_{Gr(G_x,n)})}{c(g^{*}T_{L_C^{(1^n)}(\mathcal{G})})c(j^{*}T_{L_{C'}^{(1^n)}(s)})}
\end{equation}

\begin{remark}
For any pair of flat family of sheaves $F_1$ and $F_2$ on $\mathcal{M}^{r,c}_{X, D}(n+s)\times \mathbb{C}\mathbb{P}^2$, we denote by $\mathcal{E}xt^i_p(F_1, F_2)$ the $i$th right derived functor of $p_{*}\circ \mathcal{H}om$, where $p$ is the projection from $\mathcal{M}^{r,c}_{X, D}(n+s)\times \mathbb{C}\mathbb{P}^2$ to $\mathcal{M}^{r,c}_{X, D}(n+s)$. We denote by $q$ the projection from $\mathcal{M}^{r,c}_{X, D}(n+s)\times \mathbb{C}\mathbb{P}^2$ to the other factor. If for all $x\in \mathcal{M}^{r,c}_{X, D}(n+s)$ the global Ext group $\mathrm{E}xt^i(p^{-1}(x);F_1, F_2)$ on the fiber $p^{-1}(x)$ is of constant dimension, then one can prove that $\mathcal{E}xt^i_p(F_1, F_2)$ is a vector bundle on $\mathcal{M}^{r,c}_{X, D}(n+s)$ which has the global Ext group above as the fiber over point $x$ (see \cite{art:klei}). Similar remarks apply to any closed subspace of $\mathcal{M}^{r,c}_{X, D}(n+s)$.
\end{remark}

First we prove the following technical lemma:
\begin{lem}\label{lem:tec}
Let $(F,\phi^F)$ be a framed sheaf on $X.$ Then $\mathrm{Ext}^j(F^{\vee\vee}, F(-D))=0$ for $j=0,2.$
\end{lem}
\proof
Consider the short exact sequence $0\rightarrow F \rightarrow F^{\vee\vee} \rightarrow A \rightarrow 0$ and apply the functor $\mathrm{Hom}(F^{\vee\vee}, \cdot)$, we get
\begin{eqnarray*}
0&\rightarrow& \mathrm{Hom}(F^{\vee\vee},F(-D))\rightarrow \mathrm{Hom}(F^{\vee\vee},F^{\vee\vee}(-D)) \rightarrow \mathrm{Hom}(F^{\vee\vee},A(-D))\\ 
&\rightarrow& \mathrm{Ext}^1(F^{\vee\vee}, F(-D))\rightarrow \mathrm{Ext}^1(F^{\vee\vee}, F^{\vee\vee}(-D)) \rightarrow \mathrm{Ext}^1(F^{\vee\vee}, A(-D))\\
&\rightarrow& \mathrm{Ext}^2(F^{\vee\vee}, F(-D)) \rightarrow \mathrm{Ext}^2(F^{\vee\vee}, F^{\vee\vee}(-D))\rightarrow \mathrm{Ext}^2(F^{\vee\vee}, A(-D))\rightarrow 0
\end{eqnarray*}

Since $A$ is zero-dimensional, we have $\mathrm{Ext}^i(F^{\vee\vee}, A(-D))=0$ for $i=1,2$, hence $\mathrm{Ext}^2(F^{\vee\vee}, F(-D)) \cong \mathrm{Ext}^2(F^{\vee\vee}, F^{\vee\vee}(-D)).$ By lemma \ref{lem:vanishing_framed} we get $\mathrm{Ext}^j(F^{\vee\vee}, F^{\vee\vee}(-D))=0$ for $j=0,2$, hence we get the thesis. 
\endproof

Let $\bar{\mathcal{E}}=(\bar{E}, \phi^{\bar{E}})$ be the universal family on $\mathcal{M}^{r,c}_{X, D}(n+s)\times \mathbb{C}\mathbb{P}^2$. Recall that the tangent bundle to $\mathcal{M}^{r,c}_{X, D}(n+s)$ is $\mathcal{E}xt_p^1(\bar{E}, \bar{E}\otimes q^*O_{\mathbb{C}\mathbb{P}^2}(-D))$.

We have the following result for the sheaf $h^{*}\mathcal{E}xt_p^1(\bar{E}, \bar{E}\otimes q^*O_{\mathbb{C}\mathbb{P}^2}(-D))$.
\begin{lem}
Let \(Q\) be the universal quotient bundle on the Grassmann variety $Gr(G_x, n)$. Then the full Chern class of $h^{*}\mathcal{E}xt_p^1(\bar{E}, \bar{E}\otimes q^*O_{\mathbb{C}\mathbb{P}^2}(-D))$ is equal to $(c(Q)c(Q^{\vee}))^r$, where $c(Q)$ (resp. $c(Q^{\vee})$) is the full Chern class of $Q$ (resp. of its dual).
\end{lem}
\proof
Let $\widetilde{\mathcal{E}}=(h\times \mathrm{id}_{\mathbb{C}\mathbb{P}^2})^{*}\bar{\mathcal{E}}$; then by universal property we get
\begin{equation*}
\widetilde{\mathcal{E}}\vert_{\{A\}\times \mathbb{C}\mathbb{P}^2}\cong \mathcal{F}
\end{equation*}
where $F$ is the kernel of the composition morphism $G\rightarrow G_x\rightarrow A$ for $A\in Gr(G_x, n)$ and the framed sheaf $\mathcal{F}$ is $F$ with the framing induced by inclusion to $G$.

We denote by $I_x$ the ideal sheaf of $x$ and by $\mathbb{C}_x$ the quotient sheaf $\sfrac{O_X}{I_x}$. Recall that the inclusion morphism from a torsion free sheaf on $X$ in its double dual sheaf induces a short exact sequence of sheaves on $X$; using the universal sheaf $\widetilde{\mathcal{E}}$ one can give a `universal version' of this exact sequence and obtain an exact sequence of sheaves on $Gr(G_x, n)\times \mathbb{C}\mathbb{P}^2$:
\begin{equation}\label{eq:fonda2}
0\longrightarrow \widetilde{E} \longrightarrow q^{*} G \longrightarrow B\longrightarrow 0
\end{equation}
where $B=p^{*}Q\otimes q^{*}\mathbb{C}_x$ (we denote by $p$ and $q$ the projection morphisms from $Gr(G_x, n)\times \mathbb{C}\mathbb{P}^2$ to the first factor and the second factor respectively).

We want to compute the full Chern class of the sheaf
\begin{equation*}
h^{*}\mathcal{E}xt^1_p(\bar{E},\bar{E}\otimes q^{*}O_{\mathbb{C}\mathbb{P}^2}(-D))\cong\mathcal{E}xt^1_p(\widetilde{E},\widetilde{E}\otimes q^{*}O_{\mathbb{C}\mathbb{P}^2}(-D))
\end{equation*}

Applying the functor $\mathcal{H}om_p(\cdot, \widetilde{E}\otimes q^{*}O_{\mathbb{C}\mathbb{P}^2}(-D))$ to the exact sequence \ref{eq:fonda2}, we obtain a long exact sequence of sheaves on the Grassmann variety:
\begin{eqnarray*}
0&\longrightarrow& \mathcal{E}xt^1_{p}(B, \widetilde{E}\otimes q^{*}O_{\mathbb{C}\mathbb{P}^2}(-D))\longrightarrow \mathcal{E}xt^1_{p}(q^{*} G, \widetilde{E}\otimes q^{*}O_{\mathbb{C}\mathbb{P}^2}(-D))\longrightarrow \\ 
&\longrightarrow&\mathcal{E}xt^1_{p}(\widetilde{E}, \widetilde{E}\otimes q^{*}O_{\mathbb{C}\mathbb{P}^2}(-D))\longrightarrow\mathcal{E}xt^2_{p}(B, \widetilde{E}\otimes q^{*}O_{\mathbb{C}\mathbb{P}^2}(-D))\longrightarrow 0
\end{eqnarray*}
Indeed, the sheaf $\mathcal{H}om_{p}(\widetilde{E}, \widetilde{E}\otimes q^{*}O_{\mathbb{C}\mathbb{P}^2}(-D))$ is zero since its fiber over \(A\in Gr(G_x, n)\) is $\mathrm{Hom}(F, F(-D))$, that vanishes by lemma \ref{lem:vanishing_framed}. The sheaf $\mathcal{E}xt^2_{p}(q^{*} G, \widetilde{E}\otimes q^{*}O_{\mathbb{C}\mathbb{P}^2}(-D))$ is zero because its fiber $\mathrm{Ext}^2(G, F(-D))$ over any point $A$ of the Grassmann variety, that it is zero by lemma \ref{lem:tec}.

Thus the full Chern class of $\mathcal{E}xt^1_{p}(\widetilde{E}, \widetilde{E}\otimes q^{*}O_{\mathbb{C}\mathbb{P}^2}(-D))$ is equal to:
\begin{equation*}
\frac{c(\mathcal{E}xt^2_{p}(B, \widetilde{E}\otimes q^{*}O_{\mathbb{C}\mathbb{P}^2}(-D)))c(\mathcal{E}xt^1_{p}(q^{*} G, \widetilde{E}\otimes q^{*}O_{\mathbb{C}\mathbb{P}^2}(-D)))}{c(\mathcal{E}xt^1_{p}(B, \widetilde{E}\otimes q^{*}O_{\mathbb{C}\mathbb{P}^2}(-D)))}
\end{equation*}

Now we compute the three terms in the formula above.

To compute $c(\mathcal{E}xt^1_{p}(q^{*} G, \widetilde{E}\otimes q^{*}O_{\mathbb{C}\mathbb{P}^2}(-D)))$, first we tensor the exact sequence \ref{eq:fonda2} by the invertible sheaf $q^{*} O_{\mathbb{C}\mathbb{P}^2}(-D)$ and then apply the functor $\mathcal{H}om_{p}(q^{*}G, \cdot)$. We get
\begin{eqnarray*}
0&\longrightarrow& \mathcal{H}om_{p}(q^{*}G,q^{*} G\otimes q^{*}O_{\mathbb{C}\mathbb{P}^2}(-D)) \longrightarrow \mathcal{H}om_{p}(q^{*}G,B\otimes q^{*}O_{\mathbb{C}\mathbb{P}^2}(-D))\longrightarrow\\ &\longrightarrow& \mathcal{E}xt^1_{p}(q^{*}G,\widetilde{E}\otimes q^{*}O_{\mathbb{C}\mathbb{P}^2}(-D))\longrightarrow \mathcal{E}xt^1_{p}( q^{*}G,q^{*} G\otimes q^{*}O_{\mathbb{C}\mathbb{P}^2}(-D))\longrightarrow 0
\end{eqnarray*}

The fiber over $A\in Gr(G_x,n)$ of the sheaf $\mathcal{H}om_{p}(q^{*}G,\widetilde{E}\otimes q^{*}O_{\mathbb{C}\mathbb{P}^2}(-D))$ is equal to $\mathrm{Hom}(G,F(-D))$, that vanishes by lemma \ref{lem:tec}. Thus the sheaf $\mathcal{H}om_{p}(q^{*}G,\widetilde{E}\otimes q^{*}O_{\mathbb{C}\mathbb{P}^2}(-D))$ is zero.

The sheaf $\mathcal{E}xt^1_{p}(q^{*} G, B\otimes q^{*}O_{\mathbb{C}\mathbb{P}^2}(-D)))$ is zero, because its fiber over $A\in Gr(G_x,n)$ is $\mathrm{H}^1(\mathbb{C}\mathbb{P}^2, G^{\vee}\otimes A(-D))=0$. Moreover, the sheaf $\mathcal{E}xt^1_{p}( q^{*}G,q^{*} G\otimes q^{*}O_{\mathbb{C}\mathbb{P}^2}(-D))$ is a trivial vector bundle, because its fiber over all the points of $Gr(G_{x}, n)$ is $\mathrm{Ext}^1(G,G(-D))$.

Thus
\begin{equation*}
c(\mathcal{E}xt^1_{p}(q^{*} G, \widetilde{E}\otimes q^{*}O_{\mathbb{C}\mathbb{P}^2}(-D)))=c(\mathcal{H}om_{p}(q^{*}G,B\otimes q^{*}O_{\mathbb{C}\mathbb{P}^2}(-D)))
\end{equation*}

We have
\begin{equation*}
\mathcal{H}om_{p}(q^{*}G,B\otimes q^{*}O_{\mathbb{C}\mathbb{P}^2}(-D))\cong \mathcal{H}om_{p}(q^{*}G,q^{*}(\mathbb{C}_x\otimes O_{\mathbb{C}\mathbb{P}^2}(-D)))\otimes Q
\end{equation*}
because $Q$ is a locally free sheaf.

Moreover the stalk of $\mathcal{H}om(q^{*}G,q^{*}(\mathbb{C}_x\otimes O_{\mathbb{C}\mathbb{P}^2}(-D)))$ is equal to $\mathrm{Hom}(G_x, \mathbb{C}_x)\cong \mathbb{C}^r$ over the points $(A,x)\in Gr(G_{x}, n)\times\mathbb{C}\mathbb{P}^2$ and it is equal to zero over the points $(A,y)$ with $y\neq x$. Thus the sheaf $\mathcal{H}om_{p}(q^{*}G,q^{*}(\mathbb{C}_x\otimes O_{\mathbb{C}\mathbb{P}^2}(-D)))$ is the constant sheaf $\mathbb{C}^r$ on $Gr(G_{x}, n)$.

Thus
\begin{equation*}
c(\mathcal{H}om_{p}(q^{*}G,B\otimes q^{*}O_{\mathbb{C}\mathbb{P}^2}(-D)))=c(\mathbb{C}^r\otimes Q)=c(Q)^r
\end{equation*}

Following Baranovsky's computations in the proof of lemma 6.1 in \cite{art:bara2}, we get
\begin{equation*}
c(\mathcal{E}xt^1_{p}(B, \widetilde{E}\otimes q^{*}O_{\mathbb{C}\mathbb{P}^2}(-D)))=c(\mathcal{H}om_{p}(B,B\otimes q^{*}O_{\mathbb{C}\mathbb{P}^2}(-D)))=c(Q^{\vee}\otimes Q)
\end{equation*}
and
\begin{equation*}
c(\mathcal{E}xt^2_{p}(B, \widetilde{E}\otimes q^{*}O_{\mathbb{C}\mathbb{P}^2}(-D)))=c(Q^{\vee}\otimes Q)c(Q^{\vee})^r
\end{equation*}

Summing up the results of our computation, we get
\begin{equation*}
c(\mathcal{E}xt^1_p(\widetilde{E},\widetilde{E}\otimes q^{*}O_{\mathbb{C}\mathbb{P}^2}(-D)))=\frac{c(Q^{\vee}\otimes Q)c(Q^{\vee})^rc(Q)^r}{c(Q^{\vee}\otimes Q)}=(c(Q^{\vee})c(Q))^r
\end{equation*}
\endproof

By this lemma we get
\begin{equation*}
c(h^{*}T_{\mathcal{M}^{r,c}_{X, D}(n+s)})= c(h^{*}\mathcal{E}xt_p^1(E, E\otimes q^*O_{\mathbb{C}\mathbb{P}^2}(-D)))=(c(Q^{\vee})c(Q))^r
\end{equation*}

A similar approach can be used with $T_{L_C^{(1^n)}(\mathcal{G})}$ and $T_{L_{C'}^{(1^n)}(s)}$ (as sheaves on $Gr(G_x, n)\times C$ and on $Gr(G_x, n)\times C'$). We obtain the following result:
\begin{lem}
Let $Q$ be the universal quotient bundle on the Grassmann variety $Gr(G_x, n)$. Then 
\begin{equation*}
c\left(g^{*}T_{L_C^{(1^n)}(\mathcal{G})}\right)=c\left(j^{*}T_{L_{C'}^{(1^n)}(s)}\right)=c(Q)^r
\end{equation*}
\end{lem}
\proof
Note that we can consider \(L_C^{(1^n)}(\mathcal{G})\) (resp. \(L_{C'}^{(1^n)}(s)^\bullet\)) as a variety that parametrizes Artin quotient of \(G\vert_C\) on \(C\) of length \(n\) (resp. Artin quotient of \(G'\vert_{C'}\) on \(C'\) of length \(n\), where \(\mathcal{G}'\in \mathcal{M}_{X, D}^{reg}(r,c,s)\)). From this point of view, we can use the same result of Baranovsky (see lemma 6.2 in \cite{art:bara2}) and get the thesis.
\endproof

Let $S$ be the universal subbundle on the Grassmann variety $Gr(G_x, n)$. Recall that $T_{Gr(G_x,n)}\cong S^{\vee}\otimes Q$. Now using formula \ref{eq:chern}, we get
\begin{equation*}
c(V)=\frac{(c(Q^{\vee})c(Q))^r c(S^{\vee}\otimes Q)}{c(Q)^{2r}}=\frac{c(Q^{\vee})^r c(S^{\vee}\otimes Q)}{c(Q)^{r}}
\end{equation*}

Since $c(Q)^r=c(S^{\vee}\otimes Q)c(Q^{\vee}\otimes Q)$ and $c(Q^{\vee})^r=c(S\otimes Q^{\vee})c(Q\otimes Q^{\vee})$, we obtain
\begin{equation*}
c(V)=\frac{(c(Q^{\vee}) c(S^{\vee}\otimes Q)}{c(Q)^{r}}=\frac{c(S^{\vee}\otimes Q)c(S\otimes Q^{\vee}) c(Q^{\vee}\otimes Q)}{c(S^{\vee}\otimes Q)c(Q^{\vee}\otimes Q)}=c(S\otimes Q^{\vee})
\end{equation*}

The vector bundle $S\otimes Q^{\vee}$ is the cotangent bundle on the Grassmann variety, hence we have
\begin{equation*}
(c_{n(r-n)}(V) \cap [Gr(G_x,n)])=(-1)^{n(r-n)}\int_{Gr(G_x,n)} c_{n(r-n)}(T_{Gr(G_x,n)})
\end{equation*}
and therefore
\begin{equation*}
i_{*}[L_C^{(1^n)}(\mathcal{G})] \cdot f_{*}[L_{C'}^{(1^n)}(s)] =(-1)^{n(r-n)}\chi(Gr(G_x, n))=(-1)^{(r-1)n} {r \choose n}
\end{equation*}
where we denote by $\chi(Gr(G_x, n))$ the \textit{topological Euler characteristic of} $Gr(G_x, n)$.

\subsection*{General case}

Recall that by hypothesis, $C\cap C' \cap D=\emptyset.$ Now $C$ and $C'$ intersects in $q$ points $x_1, x_2, \ldots, x_q$. Put $a_i=(-1)^{(r-1)i}{r \choose i}$.

The intersection between $L_C^{(1^n)}(\mathcal{G})$ and $L_{C'}^{(1^n)}(s)$ is, set-theoretically, the set of framed sheaves $\mathcal{F}$ such that $\mathcal{F}^{\vee\vee}\cong \mathcal{G}$ and the length of the quotient sheaf $A:=\sfrac{G}{F}$ is $n$. Note that $A$ is a quotient of $\displaystyle \bigoplus_{x\in C\cap C'} G_{x}$.

Consider a subdivision $\nu=(\nu_1, \nu_2, \ldots, \nu_q)$ of $n$, i.e. an ordered q-pla of non-negative integers $\nu_1, \ldots, \nu_q$ such that $\sum_{i=1}^{q}\nu_i=n$. Let $x=\sum_i \nu_i [x_i]$ be a cycle with $x_i\in C\cap C'$.

Let $W_\nu$ be the subset of $L_C^{(1^n)}(\mathcal{G})\cap L_{C'}^{(1^n)}(s)$ formed by all framed sheaves $\mathcal{F}$ that are kernels of morphisms $G\rightarrow A$, where $A$ is an Artin sheaf supported at the cycle $x$. Thus $W_\nu$ is isomorphic to a products of a Grassmanians $\prod_{x_i\in C\cap C'} Gr(G_{x_i}, \nu_i)$, hence it is a irreducible closed subset of $L_C^{(1^n)}(\mathcal{G})\cap L_{C'}^{(1^n)}(s)$ (note that $W_{\nu}$ does not dipend from the fixed cycle $x$ but only from the subdivision).

Moreover, by definition of $W_\nu$ follows
\begin{equation*}
L_C^{(1^n)}(\mathcal{G})\cap L_{C'}^{(1^n)}(s)=\bigcup_{\nu\: \mathrm{subdivision}}W_\nu
\end{equation*} 
and therefore the subsets $W_\nu$ are the irreducible components of $L_C^{(1^n)}(\mathcal{G})\cap L_{C'}^{(1^n)}(s)$. To each $W_\nu$ we can associate a excess bundle $V$ as before.

Since all sheaves split into direct products over the individual Grassmanians, we get that $V$ is direct sum of the excess bundles associates to each individual Grassmanian. Thus the top Chern class of $V$ is equal to the product $a_{\nu_1}\cdots a_{\nu_q}$.

By lemma 7.1 (a) and by definition 6.1.2 in \cite{book:fulton}, to get the intersection number between $L^{(1^n)}_{C}(\mathcal{G})$ and $L_{C'}^{(1^n)}(s)$ we need to sum up the products $a_{\nu_1}\cdots a_{\nu_q}$ over all subdivisions $\nu=(\nu_1, \nu_2, \ldots, \nu_q)$ of $n$. By a combinatorial argument and using the binomial formula, we get
\begin{equation*}
\sum_{n=0}^{\infty} z^{2n} <[L_{C}^{(1^n)}(\mathcal{G})],[L_{C'}^{(1^n)}(s)]>=(1+(-1)^{r-1} z^2)^{rq}.
\end{equation*}
	
\bibliography{article}

\begin{thebibliography}{10}

\bibitem{art:henni}
Amar Abdelmoubine~Henni.
\newblock Monads for torsion-free sheaves on multi-blow-ups of the projective
  plane.
\newblock 2009, \href{http://arxiv.org/abs/0903.3190}{0903.3190}.

\bibitem{art:agata}
Luis~F. Alday, Davide Gaiotto, and Yuji Tachikawa.
\newblock Liouville correlation functions from 4-dimensional gauge theories.
\newblock 2010, \href{http://arxiv.org/abs/0906.3219}{0906.3219}.

\bibitem{art:bara2}
Vladimir Baranovsky.
\newblock Moduli of sheaves on surfaces and action of the oscillator algebra.
\newblock {\em J. Differential Geom.}, 55(2):193--227, 2000.

\bibitem{art:brutanfu}
Ugo Bruzzo, Francesco Fucito, Jos{\'e}~F. Morales, and Alessandro Tanzini.
\newblock Multi-instanton calculus and equivariant cohomology.
\newblock {\em J. High Energy Phys.}, (5):054, 24 pp. (electronic), 2003.

\bibitem{art:moduli_framed}
Ugo Bruzzo and Dimitri Markushevich.
\newblock Moduli of framed sheaves on projective surfaces.
\newblock 2009, \href{http://arxiv.org/abs/0906.1436}{0906.1436}.

\bibitem{art:uhlenbeck_framed}
Ugo Bruzzo, Dimitri Markushevich, and Alexander Tikhomirov.
\newblock Uhlenbeck compactification for framed sheaves on projective surfaces.
\newblock {\em work in progress}.

\bibitem{art:brutan}
Ugo Bruzzo, Rubik Poghossian, and Alessandro Tanzini.
\newblock Poincaré polynomial of moduli spaces of framed sheaves on (stacky)
  {H}irzebruch surfaces.
\newblock 2009, \href{http://arxiv.org/abs/0909.1458}{0909.1458}.

\bibitem{art:carl}
Erik Carlsson and Andrei Okounkov.
\newblock Exts and vertex operators.
\newblock 2009, \href{http://arxiv.org/abs/0801.2565}{0801.2565}.

\bibitem{book:fulton-Young}
William Fulton.
\newblock {\em Young tableaux}, volume~35 of {\em London Mathematical Society
  Student Texts}.
\newblock Cambridge University Press, Cambridge, 1997.
\newblock With applications to representation theory and geometry.

\bibitem{book:fulton}
William Fulton.
\newblock {\em Intersection theory}, volume~2 of {\em Ergebnisse der Mathematik
  und ihrer Grenzgebiete. 3. Folge. A Series of Modern Surveys in Mathematics
  [Results in Mathematics and Related Areas. 3rd Series. A Series of Modern
  Surveys in Mathematics]}.
\newblock Springer-Verlag, Berlin, second edition, 1998.

\bibitem{art:gott2}
Lothar G{\"o}ttsche and Wolfgang Soergel.
\newblock Perverse sheaves and the cohomology of {H}ilbert schemes of smooth
  algebraic surfaces.
\newblock {\em Math. Ann.}, 296(2):235--245, 1993.

\bibitem{book:egaIV}
Alexander Grothendieck.
\newblock \'{E}l\'ements de g\'eom\'etrie alg\'ebrique. {IV}. \'{E}tude locale
  des sch\'emas et des morphismes de sch\'emas {IV}.
\newblock {\em Inst. Hautes \'Etudes Sci. Publ. Math.}, (32):361, 1967.

\bibitem{book:hl}
Daniel Huybrechts and Manfred Lehn.
\newblock {\em The geometry of moduli spaces of sheaves}.
\newblock Aspects of Mathematics, E31. Friedr. Vieweg \& Sohn, Braunschweig,
  1997.

\bibitem{phd:king}
Alastair King.
\newblock {\em Instantons and holomorphic bundles on the blown-up plane}.
\newblock PhD thesis, Worcester College (Oxford), 1989.

\bibitem{art:klei}
Steven~L. Kleiman.
\newblock Relative duality for quasicoherent sheaves.
\newblock {\em Compositio Math.}, 41(1):39--60, 1980.

\bibitem{book:mac}
Ian~G. Macdonald.
\newblock {\em Symmetric functions and {H}all polynomials}.
\newblock Oxford Mathematical Monographs. The Clarendon Press Oxford University
  Press, New York, second edition, 1995.
\newblock With contributions by A. Zelevinsky, Oxford Science Publications.

\bibitem{book:naka}
Hiraku Nakajima.
\newblock {\em Lectures on {H}ilbert schemes of points on surfaces}, volume~18
  of {\em University Lecture Series}.
\newblock American Mathematical Society, Providence, RI, 1999.

\bibitem{art:nakayo}
Hiraku Nakajima and K{\=o}ta Yoshioka.
\newblock Instanton counting on blowup. {I}. 4-dimensional pure gauge theory.
\newblock {\em Invent. Math.}, 162(2):313--355, 2005.

\bibitem{art:nekra}
Nikita~A. Nekrasov.
\newblock Seiberg-{W}itten prepotential from instanton counting.
\newblock {\em Adv. Theor. Math. Phys.}, 7(5):831--864, 2003.

\bibitem{phd:rava}
Claudio Rava.
\newblock {\em {ADHM} data for framed sheaves on {H}irzebruch surfaces}.
\newblock PhD thesis, Sissa (Trieste), 2010.

\bibitem{art:vawi}
Cumrun Vafa and Edward Witten.
\newblock A strong coupling test of {$S$}-duality.
\newblock {\em Nuclear Phys. B}, 431(1-2):3--77, 1994.

\end{thebibliography}
\bibliographystyle{hplain-hein}

\end{document}